\documentclass{extarticle}
\usepackage[utf8]{inputenc}
\usepackage{amsmath}
\usepackage{amssymb}
\usepackage{stmaryrd}
\usepackage{amsthm}
\usepackage{mathtools}
\usepackage{mathrsfs}
\usepackage{cool} 

\usepackage[nottoc]{tocbibind}
\usepackage{wasysym} 
\usepackage{listings} 
\usepackage[margin=1in]{geometry} 
\usepackage{enumitem} 
\usepackage{lastpage} 
\usepackage{hyperref} 
\usepackage{stmaryrd}

\usepackage{graphicx} 
\usepackage{tikz-cd}
\usepackage{amsmath}
\usepackage[utf8]{inputenc}
\usepackage{caption}
\usepackage{subcaption}
\usepackage[english]{babel}
\usepackage{tikz}
\usepackage{amsfonts}
\usepackage{bm}
\usepackage{bbm}
\usepackage{amssymb}
\usepackage{float}

\usepackage{fancyhdr} 
\mathchardef\mhyphen="2D 

\theoremstyle{definition}
\newtheorem{definition}{Definition}[subsection]

\newtheorem{examp}[definition]{Example}

\newtheorem{remark}[definition]{Remark}

\newtheorem{observation}[definition]{Observation}

\theoremstyle{plain}
\newtheorem{theorem}[definition]{Theorem}
\newtheorem{lemma}[definition]{Lemma}
\newtheorem{corollary}[definition]{Corollary}
\newtheorem{prop}[definition]{Proposition}
\newtheorem{thmx}{Theorem}

\newcommand\restr[2]{{
  \left.\kern-\nulldelimiterspace 
  #1 
  \vphantom{\big|} 
  \right|_{#2} 
  }}

\title{An injective Model for Twisted Derived Categories and curved Koszul Triality}

\author{Yannick Hoyer and Kristoffer Rank Rasmussen}
\newcommand{\Addresses}{{
  \bigskip
  \footnotesize

  Yannick Hoyer, \textsc{Fachbereich Mathematik, Universität Hamburg, Bundesstraße 55, 20146 Hamburg, Germany}\par\nopagebreak
  \textit{E-mail address: }\texttt{yannick.hoyer@uni-hamburg.de}

  \medskip

  Kristoffer Rank Rasmussen, \textsc{Fachbereich Mathematik, Universität Hamburg, Bundesstraße 55, 20146 Hamburg, Germany}\par\nopagebreak
  \textit{E-mail address: } \texttt{kristoffer.rasmussen@uni-hamburg.de}

  \medskip
}}

\begin{document}

\date{}
\maketitle
\begin{abstract}
Given a curved differential graded algebra $A$, we define a new model structure on the category of curved differential graded $A$-modules, called the injective Guan-Lazarev model structure. We prove that the category of CDG $A$-modules with this model structure is Quillen equivalent to the category of curved differential graded contramodules over the extended bar construction of $A$, equipped with the contraderived model structure. This result can be seen as bridging the gap between Positselski's theory of conilpotent Koszul triality and Guan-Lazarev's work on non-conilpotent Koszul duality. As an application, we use the injective Guan-Lazarev model structure to show that the tensor product is a Quillen bifunctor with respect to these model structures of the second kind. 
\end{abstract}

\maketitle
\tableofcontents

\section{Introduction}

Koszul duality for differential graded algebras (DG algebras) is a fundamental phenomenon in homological algebra and representation theory dating back to the work of Beilinson, Ginzburg and Soergel \cite{Beilinson96}. Originally thought of as a relationship between graded algebras inducing an equivalence between bounded (derived) module categories, it was observed by Positselski \cite{Positselski2011, Positselski2023KoszulDuality} that Koszul duality can be interpreted more naturally as a relationship between DG algebras and curved differential graded coalgebras (CDG coalgebras). In this more modern setting, Koszul duality identifies the derived category of a DG algebra $A$ with a certain exotic derived category of CDG comodules over a CDG coalgebra $BA$, called the bar construction of $A$. This exotic derived category, called the coderived category of CDG $BA$-comodules, arises as the homotopy category of a certain model category structure on the category of CDG $BA$-comodules, and is denoted $D^\text{co}(BA\mhyphen comod)$. We will refer to this model category structure as the coderived model structure. This approach allows for quite general versions of Koszul duality, however,  two restrictions remain. Firstly, on the algebra side, curvature cannot be permitted or else the derived category of $A$ is not defined. Secondly, $BA$ will always be a conilpotent CDG coalgebra. The first point already indicates that, in order to lift these restrictions, it is necessary to consider a different kind of derived category on the algebra side, which does not depend on the notion of cohomology. 

To solve this problem, Guan and Lazarev introduce a model structure $A\mhyphen mod_{proj}^\text{II}$ on the category $A\mhyphen mod$ of CDG $A$-modules over any CDG algebra $A$, whose homotopy category is denoted $D^\text{II}_c(A\mhyphen mod)$ (see \cite[Theorem 3.7]{GuanLazarev}) and which can be seen as an analogue to the usual projective model structure on the category of DG modules over a DG algebra. It is constructed by first considering a larger, possibly non-conilpotent, version of the classical bar construction, called the extended bar construction of $A$, denoted $\check{B}A$. There is a natural homogeneous $k$-linear map $\tau:\check{B}A\rightarrow A$ of degree 1 corresponding to a Maurer-Cartan element in the convolution algebra $Hom(\check{B}A,A)$. This gives rise to the adjunction
\begin{displaymath}
	A\otimes^\tau -:\check{B}A\mhyphen comod\rightleftarrows A\mhyphen mod:\check{B}A\otimes^\tau -,
\end{displaymath}
where the twisted tensor product $\otimes^\tau$ is the usual tensor product but equipped with a differential induced by $\tau$.

Then the model structure on the category of CDG $A$-modules is obtained by transferring the coderived model structure on the category of CDG comodules over the extended bar construction along the functor $\check{B}A\otimes^\tau -$. This transferred model structure automatically gives them a Quillen adjunction to the coderived model structure on the category of CDG $\check{B}A$-comodules, which is then shown to be a Quillen equivalence. This proves an instance of Koszul duality for CDG algebras and non-conilpotent CDG coalgebras.

It is well known that, for any CDG coalgebra $C$, one can consider the corresponding category of (left) CDG $C$-comodules. However, there is also the lesser known notion of (left) curved differential graded contramodule over $C$ (CDG $C$-contramodule), which is obtained by inverting the arrows in the characterization of CDG modules via hom spaces. Contramodules were studied extensively by Positselski \cite{Positselski2011, Positselski_2022}. Just like in the case of CDG comodules, the category of CDG contramodules over $C$ can be endowed with a model structure, whose homotopy category is called the contraderived category of CDG $C$-contramodules (see \cite[Section 8.2]{Positselski2011}), and which is denoted $D^\text{ctr}(C\mhyphen contra)$. We will refer to this model category structure as the contraderived model structure. It was shown by Positselski in \cite[Section 5.2]{Positselski2011} that there is a natural Quillen equivalence between the coderived model structure and the contraderived model structure for any CDG coalgebra $C$. This phenomenon is known as comodule-contramodule correspondence. Moreover, as was proven by Positselski in \cite[Section 6.3]{Positselski2011}, for any DG algebra $A$, a contraderived version of Koszul duality can be established, identifying the derived category of $A$ with the contraderived category of CDG $BA$-contramodules. It is also proven in \cite[Section 6.5]{Positselski2011} that both Koszul duality and comodule-contramodule correspondence are compatible with one another in the sense that, for any DG algebra $A$, there is a commutative triangle of equivalences of triangulated categories
\begin{figure}[H]
	\centering
	\begin{tikzcd}
		& & D^{co}(BA\mhyphen comod)\arrow[dd, equal]\\
		D(A\mhyphen mod)\arrow[urr, equal]\arrow[drr, equal] & &\\
		& & D^{ctr}(BA\mhyphen contra)\ .
	\end{tikzcd}
\end{figure}
This compatibility is commonly referred to as Koszul triality. It is one of the main goals of this paper to establish a curved, non-conilpotent version of Koszul triality in the framework established by Guan and Lazarev.

Let us point to our main results. Replacing CDG comodules by CDG contramodules in the setting of \cite{GuanLazarev}, the natural map $\tau:\check{B}A\to A$ gives rise to an adjunction 
\begin{displaymath}
    Hom^{\tau}(\check{B}A,-):A\mhyphen mod\rightleftarrows \check{B}A\mhyphen contra:Hom^{\tau}(A,-)\ .
\end{displaymath}
We prove that the model structure on $\check{B}A\mhyphen contra$ can be transferred to $A\mhyphen mod$ along the left adjoint $Hom^{\tau}(\check{B}A,-)$, giving rise to a model structure on the category of CDG $A$-modules. This new model structure can be seen as an analogue to the usual injective model structure on the category of DG modules over a DG algebra. 

\begin{thmx}[Theorem~\ref{theorem}, Theorem~\ref{thm:ctreq}, Theorem~\ref{thm:ctreq2}]
    Let $A$ be a CDG algebra over a field $k$ and denote by $\check{B}A$ the extended bar construction of $A$. There exists a cofibrantly generated model category structure $A\mhyphen mod_{ctr}^\text{II}$ on the category of CDG $A$-modules with closed morphisms between them, satisfying the following conditions:
    
    \begin{enumerate}
        \item An $A$-module homomorphism $f:M\rightarrow N$ is a weak equivalence if and only if
        $$f_*:Hom_A(T,M)\rightarrow Hom_A(T,N)$$
        is a quasi-isomorphism for all $T\in Tw(A)$.
        
        \item The class of cofibrations is given by the class of injections of CDG $A$-modules. In particular, every CDG $A$-module is cofibrant with respect to this model structure.
        
        \item The class of fibrations is given by the maps having the right lifting property with respect to trivial cofibrations.
    \end{enumerate}
    It is Quillen equivalent to the contraderived model structure on the category of CDG $\check{B}A$-contramodules and its homotopy category is naturally equivalent to the twisted derived category $D^\text{II}_c(A\mhyphen mod)$ defined in \cite[Theorem 3.7]{GuanLazarev}.
\end{thmx}

The model category defined above will be denoted $A\mhyphen mod_{inj}^\text{II}$.
This new model structure is then used to prove a curved, non-conilpotent version of the above mentioned Koszul triality.

\begin{thmx}[Corollary~\ref{cor:curvedtriality}]
    We have the following commutative triangle of equivalences of triangulated categories:
    \begin{figure}[H]
		\centering
		\begin{tikzcd}
			& & D^{co}(\check{B}A\mhyphen comod)\arrow[dd, equal]\\
			D_c^\text{II}(A\mhyphen mod)\arrow[urr, equal]\arrow[drr, equal] & &\\
			& & D^{ctr}(\check{B}A\mhyphen contra)\ .
		\end{tikzcd}
	\end{figure}
\end{thmx}

As an application, we establish a tensor-hom adjunction for the twisted derived category $D_c^\text{II}(A\mhyphen mod)$, which does not work without a second model structure. For any two DG algebras $A$ and $B$, we tacitly use the identification $A\mhyphen mod\mhyphen B\cong A\otimes B^{op}\mhyphen mod$, so that the above results can be applied to categories of CDG bimodules.

\begin{thmx}[Theorem~\ref{thm:tensorhom}]
    Let $A$, $B$ and $C$ be CDG algebras and consider the three model categories $A\mhyphen mod\mhyphen B_{proj}^\text{II}, B\mhyphen mod\mhyphen C_{proj}^\text{II}$ and $A\mhyphen mod\mhyphen C_{inj}^\text{II}$. The functors
    \begin{align*}
        Hom_A&:(A\mhyphen mod\mhyphen B_{proj}^\text{II})^{op}\times (A\mhyphen mod\mhyphen C_{inj}^\text{II})\rightarrow B\mhyphen mod\mhyphen C_{proj}^\text{II},\\
        Hom_C&:(B\mhyphen mod\mhyphen C_{proj}^\text{II})^{op}\times (A\mhyphen mod\mhyphen C_{inj}^\text{II})\rightarrow A\mhyphen mod\mhyphen B_{proj}^\text{II},\\
        \otimes_B&:(A\mhyphen mod\mhyphen B_{proj}^\text{II})\times (B\mhyphen mod\mhyphen C_{proj}^\text{II})\rightarrow A\mhyphen mod\mhyphen C_{inj}^\text{II}
    \end{align*}
    constitute a Quillen adjunction of two variables.
\end{thmx}

\subsection*{Acknowledgements}
We are grateful for Julian Holstein's continuous support and many helpful discussions. We would also like to thank Andrey Lazarev and Leonid Positselski for answering questions related to their work. Lastly, we appreciate the comments provided by the anonymous referee.

The first author acknowledges support by the Deutsche Forschungsgemeinschaft through SFB 1624 "Higher structures, moduli spaces and integrability" - project number 506632645 and the second author the support by Deutsche Forschungsgemeinschaft through EXC 2121 "Quantum Universe" - project number 390833306.

\section{Preliminaries}

\subsection{Notations and conventions}

Throughout this paper we will work over a fixed ground field $k$, so all vector spaces will be graded $k$-vector spaces. Unadorned tensor products and hom-spaces will always be over $k$. 

All graded objects will be cohomologically graded over a fixed abelian grading group $\Gamma$ equipped with a parity homomorphism $\Gamma\rightarrow \mathbb{Z}/2\mathbb{Z}$. Given a graded vector space $V$, we will use square brackets for shifts, i.e. $(V[n])^i=V^{n+i}$ for $i,n\in\Gamma$. If $v\in V$ is a homogeneous element, we denote by $|v|$ its degree. We will denote the linear dual of $V$ by $V^\ast$.

If $\mathcal{C}$ and $\mathcal{D}$ are categories and $F:\mathcal{C}\rightleftarrows\mathcal{D}:G$ is an adjunction, then $F$ is understood to be left adjoint to $G$.

Basic theory of model categories and triangulated categories will be assumed.

\subsection{The twisted derived category following Guan-Lazarev}

In this section, we recall the Guan-Lazarev model structure on the category of curved differential graded modules over a curved differential graded algebra and the associated homotopy category called the twisted derived category. The connection to non-conilpotent Koszul duality will be discussed briefly. 

There are no original contributions present in this section. The main reference is \cite{GuanLazarev}.

\begin{definition}
    A curved differential graded algebra (CDG algebra) $A=(A,d,h)$ is a triple consisting of the following data:
    \begin{enumerate}
        \item A graded $k$-algebra $A=\bigoplus_{i\in\Gamma}A^i$.
        
        \item A homogeneous $k$-linear map $d:A\rightarrow A$ of degree 1 satisfying the graded Leibniz rule
        \begin{displaymath}
            d(ab)=d(a)b+(-1)^{|a|}ad(b)
        \end{displaymath}
        for any homogeneous element $a\in A$ of degree $|a|$ and $b\in A$. We will refer to $d$ as the differential of $A$.
        
        \item A cocycle $h\in A^2$ satisfying $d^2(a)=ha-ah$, called the curvature element of $A$.
    \end{enumerate}

    A differential graded algebra (DG algebra) is a CDG algebra $(A,d,h)$ with $h=0$.
\end{definition}

\begin{examp}[{\cite[Page 133]{Positselski2011}}]
CDG algebras appear naturally in different contexts such as geometry, matrix factorizations and deformation theory. A CDG algebra with non-trivial curvature of geometric origin can be constructed as follows: 
    \begin{enumerate}
        \item Let $E$ be some vector bundle over a smooth affine variety $X$, together with a connection $\nabla_E:E\rightarrow E\otimes_{\mathcal{O}_X} \Omega_X^1$ i.e. a morphism of sheaves satisfying the Leibniz rule.
        \item Take the endomorphism bundle $\mathcal{E}nd(E)$ with connection $\nabla$ induced by $\nabla_E$.
        \item The algebraic de Rham algebra with coefficients in $\mathcal{E}nd(E)$, differential $\nabla$ and curvature $\nabla_E^2$ is a CDG algebra. 
    \end{enumerate}
\end{examp}

Given a graded coalgebra $C$, we denote by $\Delta$ the comultiplication and $C^\ast=Hom(C,k)$ the $k$-linear dual algebra. We will use abbreviated sweedler notation $\Delta(c)=c_{(1)}\otimes c_{(2)}$ for all $c\in C$. For any graded left $C$-comodule $M$, there is a natural left $C^\ast$-action on $N$ given by 
\begin{displaymath}
    Hom(C,k)\otimes N\to Hom(C,k)\otimes C\otimes N\to N
\end{displaymath}
induced by the left coaction on $N$ and evaluation. Similarly, any graded right $C$-comodule comes equipped with a natural right $C^\ast$-action.

\begin{definition}
A curved differential graded coalgebra (CDG coalgebra) $C=(C,d,h)$ is a triple consisting of the following data:
    \begin{enumerate}
        \item A graded $k$-coalgebra $C=\bigoplus_{i\in\Gamma}C^i$.
        
        \item A homogeneous $k$-linear map $d:A\rightarrow A$ of degree 1 satisfying 
        $$\Delta(d(c))=d(c_{(1)})\otimes c_{(2)}+(-1)^{|c_{(1)}|}c_{(1)}\otimes d(c_{(2)})$$
        for all $c\in C$.
        
        \item A homogeneous element $h\in (C^\ast)^{-2}$ satisfying $d^2(c)=hc-ch$ for all $c\in C$ and $h\circ d=0$.
    \end{enumerate}

    A differential graded coalgebra (DG coalgebra) is a CDG coalgebra $(C,d,h)$ with $h=0$.
\end{definition}

Denote by $Coalg$ the category of graded coalgebras. The forgetful functor $U:Coalg\to Vect_k$ has a right adjoint $\check{T}$ constituting an adjunction
\begin{displaymath}
    U:Coalg\rightleftarrows Vect_k:\check{T}\ .
\end{displaymath}
\begin{definition}
    Given a graded vector space $V$, we refer to $\check{T}V$ as the cofree graded coalgebra cogenerated on $V$.
\end{definition}
See \cite[section 6.2]{sweedler1969hopf} for the definition of the cofree coalgebra.
\begin{remark}
    The usual tensor coalgebra $TV$ on a vector space $V$ is always a conilpotent coalgebra. This gives rise to an adjunction
    $$U:Coalg_{conil}\rightleftarrows Vect_k:T,$$
    where $Coalg_{conil}$ denotes the category of conilpotent coalgebras.
\end{remark}

\begin{definition}[{\cite[Definition~2.5]{GuanLazarev}}]
    Given a CDG algebra $A=(A,d_A,h_A)$, choose a $k$-linear retraction of the unit map $\nu:A\to k$ and set $\bar{A}:=ker(\nu)$. Let $\check{B}A$ be the CDG coalgebra whose underlying coalgebra is $\check{T}\bar{A}[1]$ with the differential induced by the multiplication of $A$, the differential $d_A$ and the curvature element $h_A$. The curvature element of $\check{B}A$ is induced by $\nu$ and independent of the choice of $\nu$. We call $\check{B}A$ the extended bar construction of $A$.

    Similarly, given a CDG coalgebra $C=(C,d_C,h_C)$, choose a $k$-linear section $\epsilon:k\rightarrow C$ of the counit $C\rightarrow k$ and set $\overline{C}=coker(\epsilon)$. Let $\Omega C$ be the CDG algebra whose underlying graded algebra is the tensor algebra of the graded vector space $\overline{C}[-1]$ with differential induced by $d_C$ and the comultiplication of $C$. The curvature element of $\Omega C$ is induced by $\epsilon$ and is independent of the choice of $\epsilon$. We call $\Omega C$ the cobar construction of $C$.
\end{definition}
\begin{remark}
If $A$ is an augmented CDG algebra then $\check{B}A$ is a DG coalgebra. Similarly, if $C$ is a coaugmented CDG coalgebra then $\Omega C$ is a DG algebra (see \cite[Theorem 4.7]{GuanLazarev}).
\end{remark}
\begin{remark}
    In \cite{GuanLazarev}, the extended bar construction is defined as a pseudocompact CDG algebra rather than a CDG coalgebra. As the category of pseudocompact CDG algebras is anti-equivalent to the category of CDG coalgebras under the assignment $C\mapsto Hom(C,k)$ where $C\in CDG\mhyphen Coalg_k$. The inverse takes a pseudocompact CDG algebra to the continuous dual, i.e. $A\mapsto Hom_{cont}(A,k)$.
\end{remark}

\begin{theorem}[{\cite[Proposition 2.6]{GuanLazarev}}]\label{Bar-cobar adjunction}
    There is an adjunction
    \begin{align*}
        \check{B}:CDG\mhyphen Alg\rightleftarrows  CDG\mhyphen Coalg:\Omega\ ,
    \end{align*}
    where CDG-Alg denotes the category of CDG algebras and CDG-Coalg is the category of CDG coalgebras.
\end{theorem}

\begin{remark}
    In \cite{BoothLazarev2023_GlobalKoszulDuality}, the categories CDG-Alg and CDG-Coalg are equipped with model structures promoting the adjunction in Theorem \ref{Bar-cobar adjunction} to a Quillen equivalence.
\end{remark}

\begin{definition}
    Let $A$ be a CDG algebra. A curved differential graded left module over $A$ (CDG $A$-module) $M$ is a graded left $A$-module, endowed with a homogeneous $k$-linear map $d_M:M\rightarrow M$ of degree 1 satisfying
    $$d_M(am)=d_A(a)m+(-1)^{|a|}ad_M(m)$$
    and
    $$d_M^2(m)=hm\ ,$$
    for any homogeneous element $a\in A$ of degree $|a|$ and $m\in M$.

    A morphism of CDG $A$-modules is simply a morphism of graded $A$-modules. CDG $A$-modules form a locally presentable abelian DG category in the sense of \cite{Positselski_2024}. The associated abelian category of CDG $A$-modules with closed morphisms between them will be denoted $A\mhyphen mod$ and the associated homotopy category will be denoted $H(A\mhyphen mod)$.
\end{definition}

Right CDG modules can be defined analogously.

\begin{remark}
    As CDG $A$-modules over a CDG algebra $A$ are equipped with differentials that do not necessarily square to zero, there is no notion of cohomology of CDG $A$-modules. The usual construction of the derived category can therefore not be applied to the category of CDG $A$-modules.
\end{remark}

\begin{definition}
    Let $A$ be a CDG algebra. A CDG $A$-module $M$ is said to be a finitely generated twisted $A$-module if its underlying graded $A$-module is isomorphic to a finitely generated free graded $A$-module. We denote by $Tw(A)\subset H(A\mhyphen mod)$ the full subcategory of finitely generated twisted $A$-modules.
\end{definition}

\begin{remark}
    In the literature, twisted modules over DG algebras are sometimes referred to as twisted complexes, see for instance \cite{Bondal_1991}. 
\end{remark}

\begin{theorem}[{\cite[Theorem 4.6]{GuanLazarev}}]\label{thm:glmodel}
    Let A be a CDG algebra. There exists a cofibrantly generated model category structure on the category of CDG $A$-modules with closed morphisms between them, satisfying the following conditions:
    \begin{enumerate}
        \item An $A$-module homomorphism $f:M\rightarrow N$ is a weak equivalence if and only if
        $$f_*:Hom_A(T,M)\rightarrow Hom_A(T,N)\ ,$$
        is a quasi-isomorphism for all $T\in Tw(A)$.
        
        \item The class of fibrations is given by the class of surjections of CDG $A$-modules. In particular, every CDG $A$-module is fibrant with respect to this model structure.
        
        \item The class of cofibrations is given by the maps having the left lifting property with respect to trivial fibrations.
    \end{enumerate} 
\end{theorem}
\begin{definition}
    The model structure constructed in \cite[Theorem 4.6]{GuanLazarev} will be referred to as the projective Guan-Lazarev (GL) model structure and the associated model category is denoted $A\mhyphen mod^\text{II}_{proj}$. The homotopy category of the projective GL model structure is called the twisted derived category and is denoted $D_c^\text{II}(A\mhyphen mod)$.
\end{definition}

\begin{remark}
    In \cite{GuanLazarev}, the twisted derived category $D_c^\text{II}(A\mhyphen mod)$ is referred to as the compactly derived category of the second kind, as it is always a compactly generated triangulated category. The full subcategory of compact objects is given by the idempotent completion of $Tw(A)$ (see \cite[Remark 3.11]{GuanLazarev}).
\end{remark}

The twisted derived category is not an invariant of quasi-isomorphisms. This is illustrated in the next two examples.

\begin{examp}
    Let $k=\mathbb{R}$, $\mathbb{R}\mathbb{P}^2$ be the real projective plane and $C^{\ast}(\mathbb{R}\mathbb{P}^2,\mathbb{R})$ the singular cochains of $\mathbb{R}\mathbb{P}^2$. It was proven in \cite{ChuangHolsteinLazarev2021} that $D_c^\text{II}(C^{\ast}(\mathbb{R}\mathbb{P}^2,\mathbb{R})\mhyphen mod)\simeq Ind(LC_{f.g.}(\mathbb{R}\mathbb{P}^2))$, where the right hand side is the ind-completion of the category of finite dimensional $\infty$-local systems. The unit $\mathbb{R}\rightarrow C^{\ast}(\mathbb{R}\mathbb{P}^2,\mathbb{R})$ is a quasi-isomorphism, however
    $D_c^\text{II}(Vect_\mathbb{R})\simeq D(Vect_\mathbb{R})\not\simeq Ind(LC_{f.g.}(\mathbb{R}\mathbb{P}^2))$
    which shows that the unit does not induce an equivalence on the level of twisted derived categories. 
\end{examp}

\begin{examp}
    Let $A=(k[x],d)$ with $|x|=1$ and $d(x)=-x^2$. The compact objects of $D_c^\text{II}(A\mhyphen mod)$ are generated under shifts, cones and homotopy summands by $A$ and $A_x=(k[x],d_x)$ with $d_x(1)=x$. As every map from $A$ to $A_x$ and vice versa is null-homotopic, $D_c^\text{II}(A\mhyphen mod)$ splits into a product $\mathcal{S}\times \mathcal{V}$, where $\mathcal{S}$ and $\mathcal{V}$ are generated under cones, shifts and homotopy summands by $A_x$ and $A$, respectively. Clearly, the augmentation map $A\rightarrow k$ is a quasi-isomorphism, but $D_c^\text{II}(Vect_k)\simeq D(Vect_k)\not\simeq D_c^\text{II}(A\mhyphen mod)$.
\end{examp}
\begin{remark}
    The twisted derived category is an invariant of Maurer-Cartan equivalences of CDG algebras. See \cite[Definition 9.1]{BoothLazarev2023_GlobalKoszulDuality} for the definition.
\end{remark}
\begin{definition}
    Let $C=(C,d,h)$ be a CDG coalgebra. A curved differential graded left $C$-comodule (CDG $C$-comodule) $N$ is a graded $C$-comodule, endowed with a homogeneous $k$-linear map $d_N:N\rightarrow N$ of degree 1 satisfying
    \begin{align*}
        \Delta (d_N(n))=d(n_{(-1)})\otimes n_{(0)}+n_{(-1)}\otimes d_N(n_{(0)})
    \end{align*}
    and
    \begin{align*}
        d_N^2(n)=h(n_{(-1)})\otimes n_{(0)}
    \end{align*}
    for all $n\in N$.

    CDG $C$-comodules form a locally presentable abelian DG category in the sense of \cite{Positselski_2024}, which will be denoted $\text{DG}(C\mhyphen comod)$. The associated abelian category of CDG $C$-comodules with closed morphisms between them will be denoted $C\mhyphen comod$ and the associated homotopy category will be denoted $H(C\mhyphen comod)$.
\end{definition}

\begin{definition}
    Let $C$ be a CDG coalgebra. A CDG $C$-comodule $N$ is called coacyclic, if it belongs to the smallest localizing subcategory of $H(C\mhyphen comod)$ containing all totalizations of short exact sequences of CDG $C$-comodules.
\end{definition}

\begin{theorem}[{\cite[Theorem 8.2]{Positselski2011}}]\label{thm:comod}
    Let $C$ be a CDG coalgebra. There exists a cofibrantly generated model category structure on the category of CDG $C$-comodules with closed morphisms between them, satisfying the following conditions:
    \begin{enumerate}
        \item A morphism is a weak equivalence if and only if it has a coacyclic mapping cone.
    
        \item The class of cofibrations is given by the class of injections of CDG $C$-comodules.
        
        \item The class of fibrations is given by the class of surjections of CDG $C$-comodules whose kernel is a graded injective CDG $C$-comodule.
    \end{enumerate}
    The homotopy category associated to this model category structure is called the coderived category of CDG $C$-comodules and will be denoted $D^{co}(C\mhyphen comod)$.
\end{theorem}

From now on, it will always be assumed that the abelian category $C\mhyphen comod$ is equipped with the model category structure described in Theorem~\ref{thm:comod}.

The counit $\tau:\Omega\check{B}A\rightarrow A$ of the adjunction from Theorem \ref{Bar-cobar adjunction} corresponds to a Maurer-Cartan element in the convolution algebra $Hom(\check{B}A,A)$ (see \cite[Definition 3.3]{GuanLazarev}). It gives rise to a twisted tensor product functor
$$\check{B}A\otimes^\tau-:A\mhyphen mod\rightarrow \check{B}A\mhyphen comod\ ,$$
defined by $M\mapsto \check{B}A\otimes M$ with the differential $d^\tau$ defined by
\begin{align*}
    d^\tau (a\otimes b)=d_{M\otimes \check{B}A}(a\otimes b)+(m\otimes id_{\check{B}A})(id_M\otimes \tau\otimes id_{\check{B}A})(id_M\otimes\Delta)(a\otimes b)
\end{align*}
for all $a\otimes b\in M\otimes \check{B}A$. Similarly, one can define a twisted tensor product functor in the opposite direction
$$A\otimes^\tau -:\check{B}A\mhyphen comod\rightarrow A\mhyphen mod.$$

\begin{theorem}[{\cite[Theorem 4.7]{GuanLazarev}}]\label{Koszul Adjunction}
    Let $A$ be a CDG algebra. The pair of functors 
    $$A\otimes^\tau -:\check{B}A\mhyphen comod\rightleftarrows A\mhyphen mod_{proj}^\text{II}:\check{B}A\otimes^\tau -$$
    defines a Quillen equivalence. In particular, it induces an equivalence on the level of homotopy categories $D^{co}(\check{B}A\mhyphen comod)\simeq D_c^\text{II}(A\mhyphen mod)$.
\end{theorem}

\subsection{Koszul Triality}

This section is meant to collect the basic notions needed to formulate Koszul triality in the sense of Positselski and contains no original contributions. The standard reference for this topic is \cite{Positselski2011}. 

We start by giving a short introduction to curved differential graded contramodules. For a more thorough introduction to contramodules, we recommend the reference \cite{Positselski_2022}.

\begin{definition}
    Let $C$ be a graded coalgebra over $k$. A graded (left) contramodule $P$ over $C$ is a graded $k$-vector space $P=\bigoplus_{i\in\Gamma}P^i$, endowed with a homogeneous $k$-linear map $\alpha:Hom(C,P)\to P$ of degree $0$, satisfying the following contraassociativity and counity conditions:
    \begin{enumerate}
        \item The two maps 
            \begin{displaymath}
                Hom(C\otimes C,P)\to Hom(C,P)\xrightarrow{\alpha} P
            \end{displaymath}
        induced by the comultiplication of $C$ and 
            \begin{displaymath}
                Hom(C\otimes C,P)\cong Hom(C,Hom(C,P))\xrightarrow{\alpha^*} Hom(C,P)\xrightarrow{\alpha} P
            \end{displaymath}
        must be equal.

        \item The map 
            \begin{displaymath}
                P\cong Hom(k,P)\to Hom(C,P)\to P
            \end{displaymath}
        induced by the counit of $C$ must be equal to the identity on $P$.
    \end{enumerate}
\end{definition}

Given a graded coalgebra $C$ and a graded $C$-contramodule $P$, there is a natural left $C^\ast$-action on $P$, given by 
\begin{displaymath}
    Hom(C,k)\otimes P\to Hom(C,P)\to P\ .
\end{displaymath}

\begin{definition}
    Let $C$ be a CDG coalgebra. A curved differential graded contramodule over $C$ (CDG $C$-contramodule) $P$ is a graded $C$-contramodule, endowed with a homogeneous $k$-linear map $d:P\to P$ of degree 1, such that the contraaction map $\alpha:Hom(C,P)\to P$ commutes with the differentials on $Hom(C,P)$ and $P$, respectively, and such that $d^2(p)=hp$ for all $p\in P$. 

    CDG $C$-contramodules form a locally presentable abelian pretriangulated DG category in the sense of \cite{Positselski_2024}. The associated abelian category of CDG $C$-contramodules with closed morphisms between them will be denoted $C\mhyphen contra$ and the associated homotopy category will be denoted $H(C\mhyphen contra)$.
\end{definition}

\begin{definition}
    Let $C$ be a CDG coalgebra. A CDG $C$-contramodule $P$ is called contraacyclic, if it belongs to the smallest colocalizing subcategory of $H(C\mhyphen contra)$ containing all totalizations of short exact sequences of CDG $C$-contramodules.
\end{definition}

\begin{theorem}[{\cite[Theorem 8.2]{Positselski2011}}]\label{thm:contramod}
    Let $C$ be a CDG coalgebra. There exists a cofibrantly generated model category structure on the category of CDG $C$-contramodules with closed morphisms between them, satisfying the following conditions:
    \begin{enumerate}
        \item A morphism is a weak equivalence if and only if it has a contraacyclic mapping cone.
        
        \item The class of fibrations is given by the class of surjections of CDG $C$-contramodules.

        \item The class of cofibrations is given by the class of injections of CDG $C$-contramodules whose cokernel is a graded projective CDG $C$-contramodule.
    \end{enumerate}
    The homotopy category associated to this model category structure is called the contraderived category of CDG $C$-contramodules and will be denoted $D^{ctr}(C\mhyphen contra)$.
\end{theorem}

From now on, the abelian category $C\mhyphen contra$ will always be assumed to be equipped with the model category structure described in Theorem~\ref{thm:contramod}.

\begin{definition}
    Let $C$ be a CDG coalgebra. 
    \begin{enumerate}
        \item For any right CDG $C$-comodule $N$ and any left CDG $C$-contramodule $P$, define the contratensor product $N\odot_C P$ as the coequalizer of the two maps 
        \begin{displaymath}
            N\otimes Hom(C,P)\to N\otimes P
        \end{displaymath}
        induced by the left contraaction map associated to $P$ and 
        \begin{displaymath}
            N\otimes Hom(C,P)\to N\otimes C\otimes Hom(C,P)\xrightarrow{ev} N\otimes P
        \end{displaymath}
        induced by the right coaction map associated to $N$. Whenever $N$ has the structure of CDG $C$-bicomodule, $N\odot_C P$ is equipped with a natural structure of left CDG $C$-comodule (for more details see \cite[Section 2.2, Section 5.1]{Positselski2011}). We obtain a functor 
        \begin{displaymath}
            C\odot_C -:C\mhyphen contra\to C\mhyphen comod\ ,
        \end{displaymath}
        which will be denoted $\Phi_C$.

        \item For any left CDG $C$-comodule $N$, the DG vector space of comodule homomorphisms $Hom_C(C,N)$ has a natural structure of CDG $C$-contramodule via the contraaction map 
        \begin{displaymath}
            Hom(C,Hom_C(C,N))\to Hom(C,Hom(C,N))\cong Hom(C\otimes C,N)\to Hom(C,N)
        \end{displaymath}
        induced by the comultiplication of $C$ (see \cite[Section 2.2, Section 5.1]{Positselski2011}). We obtain a functor 
        \begin{displaymath}
            Hom_C(C,-):C\mhyphen comod\to C\mhyphen contra\ ,
        \end{displaymath}
        which will be denoted $\Psi_C$.
    \end{enumerate}
\end{definition}

The following theorem by Positselski is known as comodule-contramodule correspondence.

\begin{theorem}[{\cite[Theorem 5.2]{Positselski2011}}]\label{thm:posduality}
    Let $C$ be a CDG coalgebra. The pair of functors
    \begin{displaymath}
        \Phi_C:C\mhyphen contra\rightleftarrows C\mhyphen comod:\Psi_C
    \end{displaymath}
    defines a Quillen equivalence. In particular, it induces an equivalence on the level of homotopy categories $D^{ctr}(C\mhyphen contra)\cong D^{co}(C\mhyphen comod)$.
\end{theorem}

\begin{definition}
    Let $A=(A,d_A)$ be a DG algebra. Choose a $k$-linear retraction of the unit map $\nu:A\to k$ and set $\bar{A}:=ker(\nu)$. Denote by $BA:=T\bar{A}$ the tensor coalgebra on $\bar{A}$ with the differential induced by the multiplication map $\mu_A: A\otimes A\rightarrow A$ and the differential $d_A$. The curvature element of $BA$ is induced by $\nu$ and is independent of the choice of $\nu$. We will call $BA$ the bar construction of $A$.
\end{definition}

\begin{remark}
    If $A$ is an augmented DG algebra, then $BA$ is a DG coalgebra.
\end{remark}

\begin{theorem}[{\cite[Theorem 6.3]{Positselski2011}}]
    Let $A$ be a DG algebra and denote by $BA$ the bar construction of $A$. 
    \begin{enumerate}
        \item The functors $A\otimes^{\tau}-:BA\mhyphen comod\to A\mhyphen mod$ and $BA\otimes^{\tau}-:A\mhyphen mod\to BA\mhyphen comod$ induce an equivalence  $D(A\mhyphen mod)\cong D^{co}(BA\mhyphen comod)$.

        \item The functors $Hom^{\tau}(A,-):BA\mhyphen contra\to A\mhyphen mod$ and $Hom^{\tau}(BA,-):A\mhyphen mod\to BA\mhyphen contra$ induce an equivalence  $D(A\mhyphen mod)\cong D^{ctr}(BA\mhyphen contra)$.

        \item The above equivalences are compatible with the equivalence of Theorem~\ref{thm:posduality}, in the sense that they form a commutative triangle of equivalences of triangulated categories 
        \begin{figure}[H]
	    \centering
	    \begin{tikzcd}
		      & & D^{co}(BA\mhyphen comod)\arrow[dd, equal]\\
		      D(A\mhyphen mod)\arrow[urr, equal]\arrow[drr, equal] & &\\
		      & & D^{ctr}(BA\mhyphen contra)\ . 
	    \end{tikzcd}
        \end{figure}
    \end{enumerate}
\end{theorem}

\section{Injective Guan-Lazarev model structure}

\subsection{Construction}

Let $A$ be a CDG algebra, $\check{B}A$ its extended bar construction and $\tau:\check{B}A\to A$ the canonical twisting cochain. The following construction is due to Positselski \cite[Section 6.2]{Positselski2011}. There is an adjunction
\begin{align*}
    Hom^\tau(\check{B}A,-):A\mhyphen mod\rightleftarrows  \check{B}A\mhyphen contra:Hom^\tau(A,-)
\end{align*}
constructed as follows. For any CDG $A$-module $M$, let $Hom^\tau(\check{B}A,M)$ be the graded $\check{B}A$-contramodule $Hom(\check{B}A,M)$ with differential 
\begin{displaymath}
    d(f)(c)=d_M\circ f(c)-(-1)^{|f|}f\circ d_{\check{B}A}(c)+(-1)^{|f||c_{(1)}|}\tau(c_{(1)})f(c_{(2)})
\end{displaymath}
for any homogeneous element $f\in Hom(\check{B}A,M)$ of degree $|f|$ and $c\in\check{B}A$. Similarly, for any CDG $\check{B}A$-contramodule $P$ with contraaction map $\alpha:Hom(C,P)\to P$, let $Hom^{\tau}(A,P)$ be the graded $A$-module $Hom(A,P)$ with differential 
\begin{displaymath}
    d(f)(a)=d_P\circ f(a)-(-1)^{|f|}f\circ d_A(a)+\alpha(c\mapsto(-1)^{|f|+1+|c||a|}f(\tau(c)a))
\end{displaymath}
for any two homogeneous elements $f\in Hom(A,P)$ and $a\in A$ of degree $|f|$ and $|a|$, respectively.

\begin{theorem}\label{theorem}
    Let $A$ be a CDG algebra and denote by $\check{B}A$ its extended bar construction. There exists a cofibrantly generated model category structure on the category of CDG $A$-modules with closed morphisms between them, satisfying the following conditions:
    
    \begin{enumerate}
    \item An $A$-module homomorphism $f:M\rightarrow N$ is a weak equivalence if and only if
        $$f_*:Hom_A(T,M)\rightarrow Hom_A(T,N)$$
        is a quasi-isomorphism for all $T\in Tw(A)$.
        
        \item The class of cofibrations is given by the class of injections of CDG $A$-modules. In particular, every CDG $A$-module is cofibrant with respect to this model structure.
        
        \item The class of fibrations is given by the maps having the right lifting property with respect to trivial cofibrations.
    \end{enumerate}
    Moreover, the pair of functors 
    \begin{align*}
        Hom^\tau(\check{B}A,-):A\mhyphen mod_{inj}^\text{II}\rightleftarrows  \check{B}A\mhyphen contra:Hom^\tau(A,-)
    \end{align*}
    is a Quillen adjunction.
\end{theorem}

\begin{definition}
    The model structure constructed in Theorem \ref{theorem} will be referred to as the injective GL model structure and will be denoted $A\mhyphen mod^\text{II}_{inj}$.
\end{definition}

The proof of Theorem~\ref{theorem} will be achieved by using the following transfer theorem for model structures.

\begin{theorem}[{\cite[Proposition 2.2.1]{HessKedziorekRiehlShipley2017}}]\label{left-induced model}
    Consider an adjunction between locally presentable categories
    \begin{align*}
        L:\mathcal{C}\rightleftarrows \mathcal{M}:R
    \end{align*}
    where $\mathcal{M}$ is a cofibrantly generated model category with $\mathbf{W}$ being the class of weak equivalence, $\mathbf{F}$ the class of fibrations and $\mathbf{C}$ be the class of cofibrations. Suppose that the following conditions are satisfied:
    \begin{enumerate}
        \item For every object $X\in \mathcal{C}$ there exists a morphism $\epsilon_X:QX\rightarrow X$ such that $L(\epsilon_X)$ is a weak equivalence and $L(QX)$ is cofibrant.
        \item For each morphism $f:X\rightarrow Y$ in $\mathcal{C}$ there exists a morphism $Qf:QX\rightarrow QY$ satisfying $\epsilon_Y\circ Qf=f\circ \epsilon_X$.
        \item For every object $X\in \mathcal{C}$ there exists a factorisation
        \begin{align*}
            QX\coprod QX\xrightarrow[]{j} Cyl(QX)\xrightarrow[]{p}QX
        \end{align*}
        of the fold map $q:QX\coprod QX\rightarrow QX$ such that $L(j)$ is a cofibration and $L(p)$ is a weak equivalence.
    \end{enumerate}
    Then $\mathcal{C}$ can be equipped with a  cofibrantly generated model structure whose classes of cofibrations and weak equivalences are given by $L^{-1}(\textbf{C})$ and $L^{-1}(\textbf{W})$, respectively. This model structure will be referred to as the left induced model structure on $\mathcal{C}$.
\end{theorem}

\begin{lemma}\label{lem:cogen}
	Suppose $\mathcal{T}$ is a triangulated category which is generated, as a colocalizing subcategory of itself, by a class of objects $\mathcal{S}\subseteq\mathcal{T}$. Whenever $X\in\mathcal{T}$ is an object satisfying $Hom_{\mathcal{T}}(X,S)=0$ for all $S\in \mathcal{S}$, then $X$ must be zero. In other words, $\mathcal{T}$ is cogenerated by the objects in $\mathcal{S}$.
\end{lemma}
\begin{proof}
	Consider the full subcategory
	\begin{displaymath}
		\{X\}^{\perp}:=\{Y\in\mathcal{T}|Hom_{\mathcal{T}}(X,Y)=0\}\subseteq\mathcal{T}\ .
	\end{displaymath}
	This subcategory is colocalizing (see for instance \cite[Definition 2.5]{BalmerFavi2011}) and hence must contain the smallest colocalizing subcategory of $\mathcal{T}$ generated by $\mathcal{S}$. But by assumption, this is precisely $\mathcal{T}$. In particular, $X\in\{X\}^{\perp}$, so $X=0$.
\end{proof}

\begin{corollary}\label{cor:ctrcogen}
	Let $C$ be a CDG coalgebra. The triangulated category $D^{ctr}(C\mhyphen contra)$ is cogenerated by all finite dimensional CDG $C$-contramodules.
\end{corollary}
\begin{proof}
	This follows from Lemma~\ref{lem:cogen} and the fact that $D^{ctr}(C\mhyphen contra)$ is generated, as a colocalizing subcategory of itself, by all finite dimensional CDG $C$-contramodules, see \cite[Section 5.5]{Positselski2011}.
\end{proof}

The authors are grateful to the anonymous referee for pointing out the following lemma.

\begin{lemma}\label{equivalence of morphisms}
    Let $f:M\rightarrow N$ be a morphism of $A$-modules. Then
    \begin{displaymath}
			f_\ast:Hom_A(T,M)\rightarrow Hom_A(T,N)
		\end{displaymath}
    is a quasi-isomorphism for any finitely generated twisted modules $T$ if and only if
    \begin{displaymath}
			f^\ast:Hom_A(N,Hom^\tau(A,W))\rightarrow Hom_A(M,Hom^\tau(A,W))
		\end{displaymath}
    is a quasi-isomorphism for any finite dimensional $\check{B}A$-contramodule $W$.
\end{lemma}
\begin{proof}
    First notice that a finite dimensional $\check{B}A$-contramodule $W$ corresponds to a morphism of pseudocompact CDG algebras
    \begin{displaymath}
			\check{B}A^\ast\rightarrow End(W^\ast)\ .
		\end{displaymath}
    By applying \cite[Proposition 4.4, Proposition 2.6]{GuanLazarev}, we obtain a Maurer-Cartan element $\xi$ of $A\otimes End(W^\ast)$, which is nothing but a finitely generated twisted module $(A\otimes W^\ast,\xi)$. Here the differential $\xi$ is uniquely determined by the acyclic twisting cochain $\tau$ and the $\check{B}A$-contramodule structure on $W$.

    Under tensor-hom adjunction, $f^\ast$ can be expressed as a morphism
    \begin{displaymath}
			f^\ast:Hom(N\otimes W^\ast,k)^{\xi^\ast}\rightarrow Hom(M\otimes W^\ast,k)^{\xi^\ast}\ .
		\end{displaymath}
    Similarly, $f_\ast$ can be expressed as
    \begin{displaymath}
			f_\ast:(M\otimes W^\ast)^\xi\rightarrow (N\otimes W^\ast)^\xi\ .
		\end{displaymath}
    Tracing through the isomorphisms applied above, one realizes that $f^\ast$ is just the linear dual of $f_\ast$. The assertion follows from the exactness of $Hom(-,k)$.
\end{proof}
We will use the two characterizations of GL weak equivalences interchangeably. 
\begin{lemma}\label{lemma}
    Let $f:M\rightarrow N$ be a closed morphism of CDG $A$-modules. 
    \begin{enumerate}
        \item Then $f$ is a weak equivalence as in Theorem \ref{theorem} if and only if 
				\begin{displaymath}
					Hom^\tau(\check{B}A,f):Hom^\tau(\check{B}A,M)\rightarrow Hom^\tau(\check{B}A,N) 
				\end{displaymath}	
				is a weak equivalence in $\check{B}A\mhyphen contra$, i.e. an isomorphism in $D^\text{ctr}(\check{B}A\mhyphen contra)$. 
        \item Then $f$ is an injection if and only if 
				\begin{displaymath}
					Hom^\tau(\check{B}A,f):Hom^\tau(\check{B}A,M)\rightarrow Hom^\tau(\check{B}A,N) 
				\end{displaymath}
				is a cofibration in $\check{B}A\mhyphen contra$.
    \end{enumerate}
\end{lemma}
\begin{proof}
    \begin{enumerate}
        \item Let $W$ be a finite dimensional CDG $\check{B}A$-contramodule and let us denote $F=Hom^\tau(\check{B}A,-)$ and $G=Hom^\tau(A,-)$. Under the adjunction above we get the following commutative square
    $$\begin{tikzcd}
        Hom_{\check{B}A}(F(N),W)\ar[d, "\simeq"]\ar[r, "Ff^\ast"]& Hom_{\check{B}A}(F(M),W)\ar[d, "\simeq"]\\ Hom_A(N,G(W))\ar[r, "Ff_\ast"]& Hom_A(M,G(W))\ .
    \end{tikzcd}$$
    The contramodules $F(N)$ and $F(M)$ are both free as graded $\check{B}A$-contramodules and thus cofibrant, so there are natural quasi-isomorphisms of complexes of vector spaces
    \begin{displaymath}
        Hom_{\check{B}A}(F(N),W)\simeq Hom_{D^{ctr}(\check{B}A\mhyphen contra)}(F(N),W)
    \end{displaymath}
    and 
    \begin{displaymath}
        Hom_{\check{B}A}(F(M),W)\simeq Hom_{D^{ctr}(\check{B}A\mhyphen contra)}(F(M),W)\ .
    \end{displaymath}
    It follows from Corollary~\ref{cor:ctrcogen} that $F(f)$ is an isomorphism in $D^\text{ctr}(\check{B}A\mhyphen contra)$ if and only if $F(f)^\ast$ is a quasi-isomorphism for any finite dimensional CDG $\check{B}A$-contramodule. The statement now follows from the fact that $F(f)_\ast$ is a quasi-isomorphism of complexes of vector spaces precisely if $F(f)^\ast$ satisfies the same property.
    \item The functor $Hom^\tau(\check{B}A,-)$ is exact, so $Hom^\tau(\check{B}A,f)$ is an injection if and only if $f$ is an injection. In particular, if $Hom^\tau(\check{B}A,f)$ is a cofibration, then $f$ is an injection. Conversely, whenever $f$ is an injection, $Hom^\tau(\check{B}A,f)$ is an injection with cofibrant cokernel, hence a cofibration. 
    \end{enumerate}
\end{proof}
\begin{proof}[Proof of Theorem \ref{theorem}]
        The categories $A\mhyphen mod$ and $\check{B}A\mhyphen contra$ are both locally presentable \cite[Section 8.3]{Positselski2023KoszulDuality} and the model structure on $\check{B}A\mhyphen contra$ is cofibrantly generated by Theorem \ref{thm:contramod}. We verify that the three criteria of Theorem~\ref{left-induced model} are satisfied for the adjunction
        \begin{align*}
            Hom^\tau(\check{B}A,-):A\mhyphen mod\rightleftarrows  \check{B}A\mhyphen contra:Hom^\tau(A,-)\ .
        \end{align*}
        \begin{enumerate}
            \item For any CDG $A$-module $M$, set $QM=M$ and $\epsilon_M=id_M$. Then $Hom^\tau(\check{B}A,\epsilon_M)$ is a weak equivalence and $Hom^\tau(\check{B}A,M)$ is a graded free $\check{B}A$-contramodule, hence cofibrant.
            
            \item For any closed morphism of CDG $A$-modules $f:M\rightarrow N$, set $Qf=f$. Clearly $\epsilon_N\circ Qf=f\circ \epsilon_M$.
            
            \item Set $Cyl(X)=X\oplus X\oplus X[1]$ with differential 
						\begin{displaymath}
							d(a,b,c)=(d(a)+c,d(b)-c,-d(c)).
						\end{displaymath}
						There is a well-known factorization of the fold map
            \begin{displaymath}
							X\oplus X\xrightarrow{j} Cyl(X)\xrightarrow{p} X
						\end{displaymath}
            with $j$ being the inclusion $(a,b)\mapsto (a,b,0)$ and $p$ the projection $(a,b,c)\mapsto a+b$.

            Clearly, $Hom^\tau(BA,j)$ is an injection. The functor $Hom^\tau(\check{B}A,-)$ is a left adjoint, so 
            \begin{displaymath}
							coker(Hom^\tau(\check{B}A,j))\cong Hom^\tau(\check{B}A,coker(j)),
						\end{displaymath}
            which is graded free. We conclude that $Hom^\tau(\check{B}A,j)$ is a cofibration in $C\mhyphen contra$.
            According to Lemma \ref{lemma} it suffices to show that
            \begin{displaymath}
							p^\ast:Hom_A(X,Hom^\tau(A,W))\rightarrow Hom_A(Cyl(X),Hom^\tau(A,W)),
						\end{displaymath}
            is a quasi-isomorphism for all finite dimensional $\check{B}A$-contramodules $W$. As the map $p$ is a homotopy equivalence, the induced map $p^\ast$ is a quasi-isomorphism for all $W$. 
        \end{enumerate}
        It follows from Theorem \ref{left-induced model} that the left-induced model structure on $A\mhyphen mod$ exists and that the adjunction $Hom^\tau(\check{B}A,-):A\mhyphen mod_{inj}^\text{II}\rightleftarrows  \check{B}A\mhyphen contra:Hom^\tau(A,-)$ is Quillen.
        
        The characterization of the cofibrations as injections is proven in Lemma \ref{lemma} (2) and in Lemma \ref{equivalence of morphisms} it is proven that weak equivalences coincide with the weak equivalences in the projective GL model structure. 
\end{proof}

\subsection{Curved Koszul Triality}
Using the injective model structure of the previous section, we will now prove curved Koszul triality. Let $A$ be a CDG algebra. For the sake of simplicity, denote by $C:=\check{B}A$ the extended bar construction of $A$. Theorem \ref{theorem} provides the Quillen adjunction 
\begin{displaymath}
    Hom^{\tau}(C,-):A\mhyphen mod_{inj}^\text{II}\rightleftarrows C\mhyphen contra:Hom^{\tau}(A,-)\ ,    
\end{displaymath}
which we will now prove to be a Quillen equivalence.
Recall from Theorem~\ref{thm:posduality} that we have the Quillen equivalence
        $$\Phi_C: C\mhyphen contra\rightleftarrows  C\mhyphen comod:\Psi_C\ .$$

\begin{lemma}\label{lem:auxeq}
	There are natural isomorphisms of CDG $C$-comodules 
	\begin{displaymath}
		C\odot_{C} Hom^{\tau}(C,M)\cong C\otimes^{\tau}M
	\end{displaymath}	
	and CDG $C$-contramodules 
	\begin{displaymath}
		Hom_C(C,C\otimes^{\tau}M)\cong Hom^{\tau}(C,M) 
	\end{displaymath}	
	for any CDG $A$-module $M$. In other words, the following diagram commutes up to natural isomorphism.
    \begin{figure}[H]
		\centering
		\begin{tikzcd}
			& & C\mhyphen comod\arrow[dd, shift left=.1cm, "\Psi_C"]\\
			A\mhyphen mod\arrow[urr, "C\otimes^{\tau}-" ]\arrow[drr, swap, "{Hom^{\tau}(C,-)}"] & &\\
			& & C\mhyphen contra\arrow[uu, shift left=.1cm,"\Phi_C"]\ .
		\end{tikzcd}
	\end{figure}
    
\end{lemma}
\begin{proof}
	For the second isomorphism, we identify $C\otimes^{\tau}M$ as a twisted object in the DG category $\mathsf{DG}(C\mhyphen comod)$,
	\begin{displaymath}
		C\otimes^{\tau}M = (C\otimes M)(t)
	\end{displaymath}
	where $t$ is the degree one endomorphism of $C\otimes M$ given by the formula 
	\begin{displaymath}
		t(c\otimes m)=(-1)^{|c_{(1)}|}c_{(1)}\otimes\tau(c_{(2)})m\ .
	\end{displaymath}
	In other words, $C\otimes^{\tau}M$ is the CDG $C$-comodule whose underlying graded $C$-comodule is $C\otimes M$ and with differential $d_{C\otimes^\tau M}:=d_{C\otimes M}+t$. Since $Hom_C(C,-)$ is a DG functor and thus preserves all twists, we have natural isomorphisms 
	\begin{displaymath}
		Hom_C(C,C\otimes^{\tau}M)\cong Hom_C(C,(C\otimes M)(t))\cong Hom_C(C,C\otimes M)(t_*)
	\end{displaymath}
	where $t_*$ denotes the degree one endomorphism of $Hom_C(C,C\otimes M)$ induced by postcomposition with $t$. 
	
	There is a natural isomorphism of CDG $C$-contramodules 
	\begin{displaymath}
		Hom_C(C,C\otimes M)\cong Hom(C,M)
	\end{displaymath}
	(see \cite[Section 2.1]{Positselski2011}), under which the twist $t_*$ on $Hom_C(C,C\otimes M)$ gets identified with the twist $t_*'$ on $Hom(C,M)$, given by the formula
	\begin{displaymath}
		g\mapsto(c\mapsto(-1)^{|g||c_{(1)}|}\tau(c_{(1)})g(c_{(2)}))\ .
	\end{displaymath}
	Now it remains to notice that $Hom(C,M)(t_*')=Hom^{\tau}(C,M)$ by definition.
	
	The first isomorphism is constructed analogously using the canonical identification $C\odot_CHom(C,M)\cong C\otimes M$ (see \cite[Section 2.2]{Positselski2011}).
\end{proof}

Our next goal is to show that the identity on $A\mhyphen mod$ induces a Quillen equivalence $A\mhyphen mod_{proj}^\text{II}\rightleftarrows A\mhyphen mod_{inj}^\text{II}$. The proof will involve the existence of the following model structures due to Becker.

\begin{theorem}[{\cite[Proposition 1.3.6]{becker12}}]
    Let $A$ be a CDG algebra.
    \begin{enumerate}
        \item There exists a unique projective abelian model structure $A\mhyphen mod_{ctr}$ on $A\mhyphen mod$, such that the class of cofibrant objects in $A\mhyphen mod_{ctr}$ is given by the class of graded projective CDG $A$-modules.

        \item There exists a unique injective abelian model structure $A\mhyphen mod_{co}$ on $A\mhyphen mod$, such that the class of fibrant objects in $A\mhyphen mod_{co}$ is given by the class of graded injective CDG $A$-modules.
    \end{enumerate}
\end{theorem}

\begin{remark}
    In \cite{Positselski2011}, co- and contraderived categories of CDG $A$-modules are defined as certain Verdier quotients of the homotopy category of CDG $A$-module. It is an open problem whether the homotopy categories of $A\mhyphen mod_{co}$ and $A\mhyphen mod_{ctr}$ are always equivalent to Positselski's co- and contraderived categories.
\end{remark}

\begin{theorem}\label{thm:ctreq}
    The identity on $A\mhyphen mod$ induces a Quillen equivalence $A\mhyphen mod_{proj}^\text{II}\rightleftarrows A\mhyphen mod_{inj}^\text{II}$.
\end{theorem}
\begin{proof}
    Since the weak equivalences in both model structures $A\mhyphen mod_{proj}^\text{II}$ and $A\mhyphen mod_{inj}^\text{II}$ coincide, it suffices to prove that the identity is left Quillen as a functor $A\mhyphen mod_{proj}^\text{II}\to A\mhyphen mod_{inj}^\text{II}$. We first show that any cofibration $A\mhyphen mod^\text{II}_{proj}$ is an injection. It follows from \cite[Proposition 1.3.8]{becker12} that any contraacyclic CDG $A$-module in the sense of Becker (see \cite[Page 16]{becker12}) is weakly trivial in $A\mhyphen mod^\text{II}_{proj}$. In particular, the class of cofibrations in the model structure constructed in \cite[Proposition 1.3.6.1]{becker12} contains the class of cofibrations of $A\mhyphen mod^\text{II}_{proj}$. Now it is clear that $A\mhyphen mod_{proj}^\text{II}\to A\mhyphen mod_{inj}^\text{II}$ preserves cofibrations and trivial cofibrations.
\end{proof}

\begin{corollary}
    The homotopy category of the injective GL model structure is equivalent to the twisted derived category $D^\text{II}_c(A\mhyphen mod)$ defined in \cite[Theorem 3.7]{GuanLazarev}.
\end{corollary}
\begin{proof}
    As was shown in Theorem~\ref{thm:ctreq}, the identity on $A\mhyphen mod$ induces a Quillen equivalence between  $A\mhyphen mod^\text{II}_{co}$ and $A\mhyphen mod^\text{II}_{inj}$, which in turn induces an equivalence on the level of homotopy categories.
\end{proof}

\begin{remark}\label{rem:twmod}
    Denote by $Tw(A)\mhyphen mod$ the category of contravariant DG functors from $Tw(A)$ to the DG category of complexes of $k$-vector spaces. It is straightforward to check that the functor $A\mhyphen mod\to Tw(A)\mhyphen mod$ sending a CDG $A$-module $M$ to the representable DG functor $Hom_A(-,M)$ is an equivalence of categories. By construction, this equivalence identifies $A\mhyphen mod^\text{II}_{proj}$ with the classical projective model structure for the derived category $D(Tw(A)\mhyphen mod)$. Analogously, $A\mhyphen mod^\text{II}_{inj}$ is identified with the classical injective model structure for $D(Tw(A)\mhyphen mod)$. In particular, $D_c^\text{II}(A\mhyphen mod)\simeq D(Tw(A)\mhyphen mod)$. This observation was made independently by Patrick Antweiler, Julian Holstein as well as Al Guan and Andrey Lazarev.
\end{remark}

\begin{theorem}\label{thm:ctreq2}
    The Quillen adjunction 
    \begin{displaymath}
        Hom^{\tau}(C,-):A\mhyphen mod_{inj}^\text{II}\rightleftarrows C\mhyphen contra:Hom^{\tau}(A,-)
    \end{displaymath}
    provided by Theorem~\ref{theorem} is a Quillen equivalence.
\end{theorem}
\begin{proof}
    By Lemma~\ref{lem:auxeq}, the functor $Hom^{\tau}(C,-)$ factors as 
    \begin{displaymath}
        A\mhyphen mod\xrightarrow{C\otimes^{\tau}-}C\mhyphen comod\xrightarrow{Hom_C(C,-)}C\mhyphen contra\ ,
    \end{displaymath}
    so $Hom^{\tau}(C,-)$ is the right adjoint of a Quillen equivalence between the model categories $C\mhyphen contra$ and $A\mhyphen mod_{proj}^\text{II}$. In particular, the right derived functor $\mathbb{R}Hom^{\tau}(C,-)$ is an equivalence of triangulated categories $D_c^\text{II}(A\mhyphen mod)\simeq D^{ctr}(C\mhyphen contra)$. Since every object in $A\mhyphen mod_{proj}^\text{II}$ is fibrant, we have $\mathbb{R}Hom^{\tau}(C,M)=Hom^{\tau}(C,M)$ for any CDG $A$-module $M$. Now it follows from Theorem~\ref{thm:ctreq}, and the fact that every object in $A\mhyphen mod_{inj}^\text{II}$ is cofibrant, that $Hom^{\tau}(C,-)$ induces an equivalence of homotopy categories when viewed as a left Quillen functor $A\mhyphen mod_{inj}^\text{II}\to C\mhyphen contra$.
\end{proof}

\begin{corollary}\label{cor:cogen2}
    The twisted derived category $D^\text{II}_c(A\mhyphen mod)$ is cogenerated by the linear duals of all finitely generated twisted $A$-modules. Moreover, any fibrant object in $A\mhyphen mod_{inj}^\text{II}$ can be obtained as a retract of a projective limit of these cogenerators.
\end{corollary}
\begin{proof}
    It was proven in Lemma \ref{lem:cogen} that the contraderived category $D^{ctr}(C\mhyphen contra)$ is cogenerated by all finite dimensional CDG $C$-contramodules. By Theorem \ref{thm:ctreq2}, the image of these objects under the derived functor $\mathbb{R}Hom^{\tau}(A,-)$ provides a class of cogenerators of $D^\text{II}_c(A\mhyphen mod)$. Now it remains to notice that $\mathbb{R}Hom^{\tau}(A,P)=Hom^{\tau}(A,P)$ for any CDG $C$-contramodule $P$ since every object in $C\mhyphen contra$ is fibrant. As was shown in the proof of Lemma \ref{equivalence of morphisms}, any such CDG $A$-module $Hom^{\tau}(A,P)$ is the linear dual of some finitely generated twisted $A$-module.

    For the second assertion let $M$ be fibrant in $A\mhyphen mod^\text{II}_{inj}$. By Theorem~\ref{thm:ctreq2}, a fibrant replacement for $M$ is given by $N:=Hom^{\tau}(A,Hom^{\tau}(C,M))$. As $Hom^{\tau}(C,M)$ is projective as graded $C$-contramodule, it follows from \cite[Lemma A.2.3]{Positselski_2010} and \cite[Page 63]{Positselski2011} that it must be a projective limit of finite dimensional CDG $C$-contramodules. Indeed, suppose $P$ is a graded projective CDG $C$-contramodule. As described in \cite[Page 63]{Positselski2011}, there exists a natural separated and complete filtration $F$ on $P$, such that $F^nP/F^{n+1}P$ is finite dimensional for all $n$. In particular, $F^0P/F^1P=P/F^1P$ is finite dimensional, hence, by induction, $P/F^nP$ is finite dimensional for all $n$. Now it follows from \cite[Lemma A.2.3]{Positselski_2010} that $P$ is isomorphic to the projective limit of the finite dimensional CDG $C$-contramodules $P/F^nP$. We conclude that $N$ must be a projective limit of cogenerators. Since any fibrant object is a retract of its fibrant replacement, the statement follows.
\end{proof}

\begin{corollary}\label{cor:fibrant}
    Any fibrant object in $A\mhyphen mod^\text{II}_{inj}$ is a graded injective CDG $A$-module.
\end{corollary}
\begin{proof}
    This follows from Corollary \ref{cor:cogen2} since any linear dual of a finitely generated twisted $A$-module is cofree as graded $A$-module and therefore graded injective.
\end{proof} 

\begin{remark}\label{rem:quillendiag}
    To summarize, we obtain a commuting square of Quillen equivalences
    \begin{figure}[H]
        \centering
        \begin{tikzcd}
			A\mhyphen mod^\text{II}_{proj}\arrow[d, shift right=1ex]\arrow[r, shift right=1ex, "\perp"] & C\mhyphen comod \arrow[l,shift right=1ex]\arrow[d,shift left=1ex, swap, "\dashv"]\\
			A\mhyphen mod^\text{II}_{inj}\arrow[u,shift right=1ex, "\dashv"]\arrow[r, shift left=1ex] & C\mhyphen contra \arrow[l, shift left=1ex, swap, "\perp"]\arrow[u,shift left=1ex]\ .
		\end{tikzcd}
    \end{figure}
    
	If $A$ is a DG algebra, we have two more model structures on $A\mhyphen mod$, namely the injective model structure $A\mhyphen mod_{inj}$ and the projective model structure $A\mhyphen mod_{proj}$, whose homotopy category is the usual derived category $D(A\mhyphen mod)$ (see \cite[Section 3]{hinich} and \cite[Theorem 8.1]{Positselski2011}). In this case, the identity on $A\mhyphen mod$ extends the above square to the following commuting diagram of Quillen adjunctions:
	\begin{figure}[H]
		\centering
		\begin{tikzcd}
			A\mhyphen mod_{proj}\arrow[d,shift right=1ex]\arrow[r, shift left=1ex] & A\mhyphen mod^\text{II}_{proj}\arrow[l, shift left=1ex, swap, "\perp"]\arrow[d, shift right=1ex]\arrow[r, shift right=1ex, "\perp"] & C\mhyphen comod \arrow[l,shift right=1ex]\arrow[d,shift left=1ex, swap, "\dashv"]\\
			A\mhyphen mod_{inj}\arrow[u,shift right=1ex, "\dashv"]\arrow[r, shift right=1ex, "\perp"] & A\mhyphen mod^\text{II}_{inj}\arrow[l,shift right=1ex]\arrow[u,shift right=1ex, "\dashv"]\arrow[r, shift left=1ex] & C\mhyphen contra \arrow[l, shift left=1ex, swap, "\perp"]\arrow[u,shift left=1ex]\ .
		\end{tikzcd}
	\end{figure}
	Here, all vertical adjunctions are Quillen equivalences and in the left hand square both horizontal adjunctions exhibt Bousfield localisations. In particular, any DG injective DG $A$-module $M$ is a fibrant object in $A\mhyphen mod^\text{II}_{inj}$. Moreover, the horizontal adjunctions in the left square are Quillen equivalences whenever one of the following conditions is satisfied:

    \begin{enumerate}
        \item $A$ is a cofibrant DG algebra.

        \item $A$ is concentrated in non-positive degrees.

        \item $A$ is concentrated in non-negative degrees and $A^1=0$. 
    \end{enumerate}

    See \cite[page 11]{GuanLazarev}. The above situation is analogous to the case of uncurved Koszul triality described by Positselski in \cite[Section 8.4]{Positselski2011} for a DG algebra $A$ and the usual bar construction.
\end{remark}

\begin{corollary}[Curved Koszul triality]\label{cor:curvedtriality}
	We have the following commutative triangle of equivalences of triangulated categories:
	\begin{figure}[H]
		\centering
		\begin{tikzcd}
			& & D^{co}(\check{B}A\mhyphen comod)\arrow[dd, equal]\\
			D_c^\text{II}(A\mhyphen mod)\arrow[urr, equal]\arrow[drr, equal] & &\\
			& & D^{ctr}(\check{B}A\mhyphen contra)\ .
		\end{tikzcd}
	\end{figure}
\end{corollary}
\begin{proof}
    This follows by deriving the first commuting square of Remark \ref{rem:quillendiag}.
\end{proof}

\subsection{Abelian model structures}

The goal of this section is to prove that both model structures $A\mhyphen mod^\text{II}_{proj}$ and $A\mhyphen mod^\text{II}_{inj}$ are abelian, thus obtaining an explicit description of the cofibrations in $A\mhyphen mod^\text{II}_{proj}$ and the fibrations in $A\mhyphen mod^\text{II}_{inj}$, respectively. 

We recall the notion of abelian model structures.
\begin{definition}
    Let $\mathcal{A}$ be an abelian category. We say that a model structure on $\mathcal{A}$ is abelian, if the following two conditions are satisfied:
    \begin{enumerate}
        \item A morphism is a cofibration if and only if it is an injection with cofibrant cokernel.

        \item A morphism is a fibration if and only if it is a surjection with fibrant kernel.
    \end{enumerate}
    An abelian model structure on $\mathcal{A}$ is called injective (projective, respectively) if every object of $\mathcal{A}$ is cofibrant (fibrant, respectively)
\end{definition}

\begin{lemma}\label{lem:abmod}
    Let $A$ be a CDG algebra.
    \begin{enumerate}
        \item Any cofibration in $A\mhyphen mod^\text{II}_{proj}$ is an injection with cofibrant cokernel.

        \item Any fibration in $A\mhyphen mod^\text{II}_{inj}$ is a surjection with fibrant kernel.
    \end{enumerate}
\end{lemma}
\begin{proof}
    We only prove the second statement. The proof for the first statement is analogous. Let $f:M\to N$ be a fibration in $A\mhyphen mod^\text{II}_{inj}$. It follows from Theorem~\ref{thm:ctreq} that any fibration in $A\mhyphen mod^\text{II}_{inj}$ is also a fibration in $A\mhyphen mod^\text{II}_{proj}$, hence a surjection.

    We claim that $ker(f)$ is fibrant. Let $p:X\to Y$ be a trivial cofibration and consider the following lifting problem.
    \begin{figure}[H]
        \centering
        \begin{tikzcd}
            X\arrow[r]\arrow[d] & ker(f)\arrow[d,"p"]\\
            Y\arrow[r] & 0\ .
        \end{tikzcd}
    \end{figure}
    Denote by $\iota:ker(f)\to M$ the canonical inclusion and consider the following commutative diagram. 
    \begin{figure}[H]
        \centering
        \begin{tikzcd}
            X\arrow[r]\arrow[d] & ker(f)\arrow[d,"p"]\arrow[r,"\iota"] & M\arrow[d,"f"]\\
            Y\arrow[r] & 0\arrow[r] & N\ .
        \end{tikzcd}
    \end{figure}
    As $f$ is a fibration, there exists a lift $h':Y\to M$, making the outer diagram commute. In particular, $f\circ h'=0$, so $h'$ factors through $ker(f)$, providing the desired lift $h:Y\to ker(f)$.
\end{proof}

To show the converse inclusions, we will use the following lifting criterion for injections and surjections in abelian categories. 

\begin{lemma}[{\cite[Lemma 9.1.1]{Positselski_2010}}]\label{lem:poslift}
    Let $\mathcal{A}$ be an abelian category, $i:X\to Y$ an injection in $\mathcal{A}$ and $p:M\to N$ a surjection in $\mathcal{A}$. If $Ext^1(coker(i),ker(p))=0$, then any lifting problem of the following form admits a solution.
    \begin{figure}[H]
        \centering
        \begin{tikzcd}
            X\arrow[r]\arrow[d,"i"] & M\arrow[d,"f"]\\
            Y\arrow[r] & N\ .
        \end{tikzcd}
    \end{figure}
\end{lemma}

\begin{lemma}\label{lem:abmod2}
    Let $M\in A\mhyphen mod$ be a CDG $A$-module whose image in $D_c^\text{II}(A\mhyphen mod)$ is zero. 
    \begin{enumerate}
        \item For any cofibrant object $X$ in $A\mhyphen mod^\text{II}_{proj}$ we have $Ext^1(X,M)=0$.

        \item For any fibrant object $Y$ in $A\mhyphen mod^\text{II}_{inj}$ we have $Ext^1(M,Y)=0$.
    \end{enumerate}
\end{lemma}
\begin{proof}
    We only prove the second statement. The first assertion is proven dually. 

    Suppose first that $Y$ is the linear dual of a finitely generated twisted $A$-module. Suppose we are given a short exact sequence of CDG $A$-modules
    \begin{displaymath}
        0\to Y\xrightarrow{f} X\xrightarrow{g} M\to 0\ .
    \end{displaymath}
    As $Y$ is graded injective, applying the DG functor $Hom_A(-,Y)$ to the above short exact sequence yields a short exact sequence of complexes of vector spaces 
    \begin{displaymath}
        0\to Hom_A(M,Y)\xrightarrow{g^\ast} Hom_A(X,Y)\xrightarrow{f^\ast} Hom_A(Y,Y)\to 0\ .
    \end{displaymath}
    By assumption, $Hom_A(M,Y)$ is acyclic, so $f^\ast$ must be a quasi-isomorphism by the long exact sequence in cohomology. In particular, there exist a closed morphism $\phi:X\to Y$ of degree $0$ as well as a morphism $\psi:Y\to Y$ of degree $-1$, such that $id_Y=\phi\circ f+d(\psi)$. As $f^\ast$ is surjective, we can choose some morphism of CDG $A$-modules $\psi':X\to Y$ satisfying $\psi'\circ f=\psi$. A straightforward calculation shows that $(\phi-d(\psi'))\circ f=id_Y$, hence the original short exact sequence splits. 

    The general case now follows from the observation that any fibrant object in the injective GL model category is a limit of a projective system of linear duals of finitely generated twisted modules with each transition map being a surjection, see Corollary~\ref{cor:cogen2}. If we replace a fibrant object $Y$ in the above argument with a projective limit as described above, then $Hom_A(M,Y)$ is a limit of a projective system of acyclic complexes of vector spaces such that each transition map is surjective, thus $Hom_A(M,Y)$ is acyclic by \cite[Theorem 3.1]{ROOS_2006}.
\end{proof}

\begin{theorem}\label{thm:abmod}
    $A\mhyphen mod^\text{II}_{proj}$ is a projective abelian model structure and $A\mhyphen mod^\text{II}_{inj}$ is an injective abelian model structure.
\end{theorem}
\begin{proof}
    We only prove the statement for $A\mhyphen mod^\text{II}_{inj}$. The other case is proven dually. 

    The fact that every object in $A\mhyphen mod^\text{II}_{inj}$ is cofibrant is clear, since every injection is a cofibration. This implies that the class of cofibrations in $A\mhyphen mod^\text{II}_{inj}$ equals the class of injections of CDG $A$-modules with cofibrant cokernel. 

    By Lemma \ref{lem:abmod}, any fibration in $A\mhyphen mod^\text{II}_{inj}$ is a surjection with fibrant kernel. Conversely, let $f:M\to N$ be a surjection of CDG $A$-modules with fibrant kernel and consider a lifting problem of the form
    \begin{displaymath}
        \begin{tikzcd}
            X\arrow[r]\arrow[d,"i"] & M\arrow[d,"f"]\\
            Y\arrow[r] & N
        \end{tikzcd}
    \end{displaymath}
    with $i:X\to Y$ being an acyclic cofibration. We claim that the image of $coker(i)$ in $D_c^\text{II}(A\mhyphen mod)$ is zero. Indeed, let $Q$ be the linear dual of a finitely generated twisted $A$-module. Applying the functor $Hom_A(-,Q)$ to the short exact sequence 
    \begin{displaymath}
        0\to X\xrightarrow{i}Y\to coker(i)\to 0
    \end{displaymath}
    yields an exact sequence of complexes of vector spaces
    \begin{displaymath}
        0\to Hom_A(coker(i),Q)\to Hom_A(X,Q)\xrightarrow{i^\ast} Hom_A(Y,Q)\to 0\ .
    \end{displaymath}
    As $i$ is a weak equivalence, $i^\ast$ is a quasi-isomorphism, so $Hom_A(coker(i),Q)$ is acyclic by the long exact sequence in cohomology. By Lemma \ref{lem:abmod2} we have $Ext^1(coker(i), Q)=0$, so the above lifting problem admits a solution by Lemma \ref{lem:poslift}.
\end{proof}

\begin{remark}
    Theorem~\ref{thm:abmod} also follows directly from Remark~\ref{rem:twmod} since the classical model structures on $Tw(A)\mhyphen mod$ are abelian. Nonetheless, the authors decided to include a more explicit computation as it adds valuable insight into the GL model structures.
\end{remark}

\begin{corollary}
    Let $A$ be a CDG algebra.
    \begin{enumerate}
        \item The identity functor $A\mhyphen mod_{ctr}\to A\mhyphen mod^\text{II}_{proj}$ defines a right Bousefield localization.

        \item The identity functor $A\mhyphen mod_{co}\to A\mhyphen mod^\text{II}_{inj}$ defines a left Bousefield localization.
    \end{enumerate}
\end{corollary}
\begin{proof}
    By \cite[Proposition 4.10(1)]{GuanHolsteinLazarev2025}, any cofibrant object in $A\mhyphen mod^\text{II}_{proj}$ is graded projective, hence also cofibrant in $A\mhyphen mod_{ctr}$. 

    Dually, by Corollary~\ref{cor:fibrant}, any fibrant object in $A\mhyphen mod^\text{II}_{inj}$ is graded injective. 
\end{proof}

\subsection{Derived GL Tensor-Hom Adjunction}

As an application of the injective GL model structure, we establish a tensor-hom adjunction for the twisted derived category and provide an example that illustrates the necessity of two distinct model structures for this to work. Moreover, we prove that curved Koszul triality is compatible with various binary operations on CDG modules, CDG comodules and CDG contramodules.

We recall the notion of an adjunction of two variables.

\begin{definition}[\cite{Hovey1999}]
    Let $\mathcal{C},\mathcal{D}$ and $\mathcal{E}$ be categories. An adjunction of two variables is given by a tuple $(\boxtimes, Hom_r,Hom_l,\varphi_r,\varphi_l)$ where $\boxtimes:\mathcal{C}\times \mathcal{D}\rightarrow \mathcal{E}$, $Hom_r:\mathcal{D}^{op}\times \mathcal{E}\rightarrow \mathcal{C}$ and $Hom_r:\mathcal{C}^{op}\times \mathcal{E}\rightarrow \mathcal{D}$ are functors and $\varphi_r$ and $\varphi_l$ are natural isomorphisms
    \begin{align*}
        Hom_\mathcal{C}(X,Hom_r(Y,Z))\xleftarrow{\varphi_r} Hom_\mathcal{E}(X\boxtimes Y,Z)\xrightarrow{\varphi_l} Hom_\mathcal{D}(Y,Hom_l(X,Z))
    \end{align*}
    for $X\in\mathcal{C}$, $Y\in\mathcal{D}$ and $Z\in\mathcal{E}$.

    The natural isomorphisms $\phi_r$ and $\phi_l$ will not be mentioned explicitly if they are clear from context.
\end{definition}

\begin{definition}
    Let $A$ and $B$ be CDG algebras. A CDG $A\mhyphen B$-bimodule is a CDG $A\otimes B^{op}$-module. We will denote by $A\mhyphen mod\mhyphen B$ the category of CDG $A\mhyphen B$-bimodules with closed morphisms between them.
\end{definition}

In particular, for any two CDG algebras $A$ and $B$ we have the two model structures $A\mhyphen mod\mhyphen B_{proj}^\text{II}$ and $A\mhyphen mod\mhyphen B_{inj}^\text{II}$ constructed in \cite[Theorem 4.6]{GuanLazarev} and \ref{theorem}, respectively.

\begin{examp}
    Let $A$, $B$ and $D$ be CDG algebras, $M$ an $A\mhyphen B$-bimodule and $N$ a $B\mhyphen D$-bimodule. The tensor product $M\otimes_B N$ is defined analogously to how it is defined for DG modules.
    The tensor product $\otimes_B$ together with the two functors
    \begin{align*}
        Hom_A&:(A\mhyphen mod\mhyphen B)^{op}\times (A\mhyphen mod\mhyphen D)\rightarrow B\mhyphen mod\mhyphen D\\
        Hom_D&:(B\mhyphen mod\mhyphen D)^{op}\times (A\mhyphen mod\mhyphen D)\rightarrow A\mhyphen mod\mhyphen B
    \end{align*}
    constitute an adjunction of two variables. See \cite[Section 3.10]{Positselski2011} for an in-depth discussion of the tensor product of CDG modules.
\end{examp}

\begin{definition}\label{Quillen Bifunctors}
    Let $\mathcal{C},\mathcal{D}$ and $\mathcal{E}$ be model categories. An adjunction of two variables $(\boxtimes, Hom_r,Hom_l)$ is a Quillen adjunction of two variables if the following two conditions are satisfied.
    \begin{enumerate}
        \item Given cofibrations $f:U\rightarrow V$ in $\mathcal{C}$ and $g:W\rightarrow X$ in $\mathcal{D}$, the pushout product
        \begin{align*}
            f\ \square\  g:(V\boxtimes W)\coprod_{U\boxtimes W}(U\boxtimes X)\rightarrow V\boxtimes X
        \end{align*}
        is a cofibration in $\mathcal{E}$.
        \item Given cofibrations $f$ in $\mathcal{C}$ and $g$ in $\mathcal{D}$, if either $f$ or $g$ is a trivial cofibration in $\mathcal{C}$ or $\mathcal{D}$, respectively, then $f\ \square\ g$ is a trivial cofibration in $\mathcal{E}$. 
    \end{enumerate}
\end{definition}

\begin{lemma}[{\cite[Theorem 4.6]{GuanLazarev}}]\label{trivial cof.}
    Let $A$ and $B$ be CDG algebras. The projective GL model structure on $A\mhyphen mod\mhyphen B$ is cofibrantly generated by the following classes of maps.
    \begin{enumerate}
        \item The class of generating cofibrations is given by maps between finitely generated twisted CDG $A\mhyphen B$-bimodules of the form
    \begin{align*}
        \text{id}_A\otimes i\otimes \text{id}_B:A\otimes V\otimes B\rightarrow A\otimes V'\otimes B
    \end{align*}
    where $i:V\rightarrow V'$ is an injection of graded vector spaces.
    \item The class of generating trivial cofibrations is given by the generating cofibrations that are also weak equivalences.
    \end{enumerate}
\end{lemma}
\begin{corollary}\label{trivial cof. are homotopy eq.}
    Let $A$ and $B$ be CDG algebras. Any generating trivial cofibration in $A\mhyphen mod\mhyphen B_{proj}^\text{II}$ is a homotopy equivalence.
\end{corollary}
\begin{proof}
     The finitely generated twisted modules are all bifibrant objects in the projective GL model structure. Weak equivalences between bifibrant objects are homotopy equivalences. 
\end{proof}

\begin{lemma}\label{lem:pushout}
    Let $A$, $B$ and $D$ be CDG algebras, $f:U\rightarrow V$ a map of CDG $A\mhyphen B$-bimodules and $g:W\rightarrow X$ a map of CDG $B\mhyphen D$-bimodules. Consider the following pushout diagram:
$$\begin{tikzcd}[sep=huge]
    U\otimes_B W\ar[r, "id_U\otimes_B g"]\ar[d, "f\otimes_B id_W"'] & U\otimes_B\ar[d, "i_2"] X\\
    V\otimes_B W\ar[r, "i_1"] & Z \ . \ar[ul, phantom, very near start, "\ulcorner"]
\end{tikzcd}$$
There is a natural isomorphism of CDG $A\mhyphen D$-bimodules
$$Z\cong (V\otimes_B W)\oplus (U\otimes_B X)/\sim$$
where $(f(u)\otimes_B w,0)\sim (0,u\otimes_B g(w))$ for all $u\in U$ and $w\in W$. 
\end{lemma}
\begin{proof}
    There are two natural inclusions $\iota_1:V\otimes_B W\rightarrow (V\otimes_B W)\oplus (U\otimes_B X)/\sim$ and $\iota_2:U\otimes_B X\rightarrow (V\otimes_B W)\oplus (U\otimes_B X)/\sim$. 

    Denote by 
    \begin{align*}
        \varphi:(V\otimes_B W)\oplus (U\otimes_B X)/\sim&\rightarrow Z\\
        (a,b)&\mapsto i_1(a)+i_2(b)
    \end{align*}
    the unique map of CDG $A\mhyphen D$-bimodules compatible with the inclusions. Invoking the universal property of $Z$, we get a map
    \begin{align*}
       \psi:Z\rightarrow (V\otimes_B W)\oplus (U\otimes_B X)/\sim\ ,
    \end{align*}
    which turns out to be the two-sided inverse $\varphi^{-1}$.
\end{proof}

\begin{observation}
    Under the identification of Lemma~\ref{lem:pushout}, the pushout product of two maps $f:U\to V$ and $g:W\to X$ is explicitly given by $(v\otimes_B w,u\otimes_B x)\mapsto v\otimes_B g(w)+f(u)\otimes_B x$ for all $v\in V$, $u\in U, w\in W$ and $x\in X$.
\end{observation}

\begin{lemma}
    Let $A$, $B$ and $D$ be CDG algebras. Let $f,g:M\rightarrow N$ be maps of $A\mhyphen B$-bimodules and assume that $f$ and $g$ are homotopic as $A\mhyphen B$-bimodule maps. Let $K$ be a $B\mhyphen D$-bimodule, then $f\otimes_B K$ and $g\otimes_B K$ are homotopic as $A\mhyphen D$-bimodule maps.
\end{lemma}
\begin{proof}
    Let $h:M\rightarrow N$ be homotopy between $f$ and $g$, i.e. $h$ a degree -1 map satisfying
    \begin{align*}
        f-g=d_N h+h d_M\ .
    \end{align*}
    It can easily be verified that $h\otimes_B \text{id}_K$ is a homotopy between $f\otimes_B \text{id}_K$ and $g\otimes_B \text{id}_K$.
\end{proof}

\begin{theorem}\label{thm:tensorhom}
    Let $A$, $B$ and $D$ be CDG algebras and consider the three model categories $A\mhyphen mod\mhyphen B_{proj}^\text{II}$, $B\mhyphen mod\mhyphen D_{proj}^\text{II}$ and $A\mhyphen mod\mhyphen D_{inj}^\text{II}$. The functors
    \begin{align*}
        Hom_A&:(A\mhyphen mod\mhyphen B_{proj}^\text{II})^{op}\times (A\mhyphen mod\mhyphen D_{inj}^\text{II})\rightarrow B\mhyphen mod\mhyphen D_{proj}^\text{II}\ ,\\
        Hom_D&:(B\mhyphen mod\mhyphen D_{proj}^\text{II})^{op}\times (A\mhyphen mod\mhyphen D_{inj}^\text{II})\rightarrow A\mhyphen mod\mhyphen B_{proj}^\text{II}\ ,\\
        \otimes_B&:(A\mhyphen mod\mhyphen B_{proj}^\text{II})\times (B\mhyphen mod\mhyphen D_{proj}^\text{II})\rightarrow A\mhyphen mod\mhyphen D_{inj}^\text{II}
    \end{align*}
    constitute a Quillen adjunction of two variables. In particular, the corresponding derived functors constitute an adjunction of two variables
    \begin{align*}
        \mathbb{R}Hom_A^\text{II}&:D_c^\text{II}(A\mhyphen mod\mhyphen B)^{op}\times D_c^\text{II}(A\mhyphen mod\mhyphen D)\rightarrow D_c^\text{II}(B\mhyphen mod\mhyphen D)\ ,\\
        \mathbb{R}Hom_D^\text{II}&:D_c^\text{II}(B\mhyphen mod\mhyphen D)^{op}\times D_c^\text{II}(A\mhyphen mod\mhyphen D)\rightarrow D_c^\text{II}(A\mhyphen mod\mhyphen B)\ ,\\
        \otimes_B^{\text{II},\mathbb{L}}&:D_c^\text{II}(A\mhyphen mod\mhyphen B)\times D_c^\text{II}(B\mhyphen mod\mhyphen D)\rightarrow D_c^\text{II}(A\mhyphen mod\mhyphen D)\ .
    \end{align*}
\end{theorem}
\begin{proof}
We have to show the two criteria listed in Definition \ref{Quillen Bifunctors} are fulfilled. Let $f:U\rightarrow V$ and $g:W\rightarrow X$ be cofibrations in $A\mhyphen mod\mhyphen B_{proj}^\text{II}$ and $B\mhyphen mod\mhyphen D_{proj}^\text{II}$, respectively. The projective GL model structure is cofibrantly generated, so according to \cite[Lemma 4.2.4]{Hovey1999}, we can assume that $f$ and $g$ are generating cofibrations. 
\begin{enumerate}
    \item 
    Finitely generated twisted $A\mhyphen B$-bimodules are free right $B$-modules and similarly finitely generated twisted $B\mhyphen C$-bimodules are free left $B$-modules, thus $f\otimes_B W$ and $U\otimes_B g$ are injections. Consider the $A\mhyphen D$-bimodule homomorphism
    \begin{align*}
        \phi:(V\otimes_B W)\oplus (U\otimes_B X)&\rightarrow V\otimes X\\
        (v\otimes_B w,u\otimes_B x)&\mapsto v\otimes_B g(w)+f(u)\otimes_B x
    \end{align*}
    and note $(V\otimes_B W\oplus U\otimes_B X)/\text{ker}(\phi)\cong (V\otimes_B W)\oplus (U\otimes_B X)/\sim$. We conclude that $f\ \square\ g$ is an injection and, in particular, a cofibration in $A\mhyphen mod\mhyphen D_{inj}^\text{II}$.
    
    \item  Assume that $f$ is a generating trivial cofibration, then, according to Lemma \ref{trivial cof.} and Corollary \ref{trivial cof. are homotopy eq.}, $f$ is a homotopy equivalence. Let $f'$ denote the homotopy inverse. We claim that $f\ \square \ g$ is a homotopy equivalence, with homotopy inverse given by the map
    \begin{align*}
        \psi:V\otimes_B X&\rightarrow (V\otimes_B W)\oplus (U\otimes_B X)/\sim\\
        v\otimes x&\mapsto (0,f'(v)\otimes_B x)\ .
    \end{align*}
    Indeed, the composition $(f\ \square \ g)\circ \psi$ is given by $v\otimes_B x\mapsto f\circ f'(v)\otimes_B x$, thus $(f\ \square \ g)\circ \psi$ is homotopic to $\text{id}_{V\otimes_B X}$. Conversely, the composition $\psi\circ (f\ \square \ g)$ is given by
    \begin{align*}
        (v\otimes_B w,u\otimes_B x)&\mapsto(0,f'(v)\otimes_B g(w)+f'\circ f(u)\otimes_B x)\\&=(f\circ f'(v)\otimes_B w,f'\circ f(u)\otimes_B x)\ ,
    \end{align*}
    thus $\psi\circ (f\ \square \ g)$ is homotopic to $\text{id}_{(V\otimes_B W)\oplus (U\otimes_B X)/\sim}$.
\end{enumerate}
\end{proof}
The injective GL model structure remedies the standing problem that 
\begin{align*}
    \otimes_B&:A\mhyphen mod\mhyphen B^\text{II}_{proj}\times B\mhyphen mod\mhyphen D^\text{II}_{proj}\rightarrow A\mhyphen mod\mhyphen D^\text{II}_{proj}
\end{align*}
is generally not a Quillen bifunctor. This is evident from the example below.
\begin{examp}\label{ex:notcofib}
    Let $k$ be an algebraically closed field, $A=k[x]$ with $|x|=0$ and $B=k\langle \epsilon\rangle$ with $|\epsilon|=1$. 
    
		Consider the $B\mhyphen A$-bimodule $X=(k\langle \epsilon\rangle\otimes k[x],d_X)$ with $d_X(1\otimes 1)=\epsilon\otimes x$. Clearly, $A$ is cofibrant in $A\mhyphen mod_{proj}^\text{II}$ and $X$ is cofibrant in $B\mhyphen mod\mhyphen A_ {proj}^\text{II}$. We want to show that $X\otimes_A A$ cannot be cofibrant in $B\mhyphen mod^\text{II}_{proj}$. Let ${_BX}$ denote $X$ restricted to its $B$-module structure and assume that $X\otimes_A A\cong{_BX}$ is cofibrant. The closed map $\psi:k[-1]\rightarrow{_BX}$ defined by $1\mapsto \epsilon\otimes 1$ has a coacyclic cone \cite[page 32]{Positselski2011} and is therefore a weak equivalence \cite[Proposition 4.10 (2)]{GuanHolsteinLazarev2025}. The trivial $B$-module $k$ can be cofibrantly replaced
    \begin{align*}
        \varphi:\bigoplus_{\lambda\in k} \left(...\xrightarrow{\epsilon}k\langle \epsilon\rangle_\lambda\xrightarrow{\epsilon} k\langle \epsilon\rangle_\lambda\xrightarrow{\epsilon}k\langle \epsilon\rangle_\lambda \right)[-1]\rightarrow k[-1]
    \end{align*}
where $k\langle \epsilon\rangle_\lambda$ is the free rank one $B$-module with differential given by multiplication by $\lambda \epsilon$. The map $\varphi$ is defined component-wise as the projection
$(...,b_3 \epsilon+a_3,b_2 \epsilon+a_2,b_1 \epsilon+a_1)\mapsto a_1,$
where $a_i,b_i\in k$. It is a straightforward computation to verify that $\varphi$ is a weak equivalence using the fact that $D^\text{II}_c(B)$ is compactly generated by the finitely generated twisted rank one modules $k\langle \epsilon\rangle_\lambda$ (see \cite[Example 3.12]{GuanLazarev}).

The map $\psi\circ \varphi$ is a weak equivalence of cofibrant objects, thus it is a homotopy equivalence. Let $\rho:{_BX}\rightarrow \bigoplus_{\lambda\in k} \left(...\xrightarrow{\epsilon}k\langle \epsilon\rangle_\lambda\xrightarrow{\epsilon} k\langle \epsilon\rangle_\lambda\xrightarrow{\epsilon}k\langle \epsilon\rangle_\lambda \right)[-1]$ be the homotopy inverse of $\psi\circ \varphi$. Every generator of $_BX$ is in degree 0, while 
\begin{displaymath}
	\bigoplus_{\lambda\in k} \left(...\xrightarrow{\epsilon}k\langle \epsilon\rangle_\lambda\xrightarrow{\epsilon} k\langle \epsilon\rangle_\lambda\xrightarrow{\epsilon}k\langle \epsilon\rangle_\lambda \right)[-1]
\end{displaymath}	
is concentrated in positive degrees, so $\rho=0$. Hence $\rho$ cannot be a homotopy equivalence. This contradicts the assumption that $_BX$ is cofibrant.
\end{examp}

\begin{remark}
    Example \ref{ex:notcofib} implies that when the target $A\mhyphen mod\mhyphen D$ of $\otimes_B$ is equipped with the projective GL model structure, then $\otimes_B$ is not generally a Quillen bifunctor. The authors do not expect 
    \begin{align*}
    \otimes_B&:A\mhyphen mod\mhyphen B^\text{II}_{inj}\times B\mhyphen mod\mhyphen D^\text{II}_{proj}\rightarrow A\mhyphen mod\mhyphen D^\text{II}_{inj}
    \end{align*}
    to be a Quillen bifunctor either, yet are unaware of a counterexample. 
\end{remark}

\begin{remark}
	Note that the derived tensor product $\otimes_B^{\text{II},\mathbb{L}}$ is not proven to be associative in any way. It is an open question whether an appropriate associativity holds.
\end{remark}

Having established a derived tensor product on twisted derived categories, we consider binary operations on co- respectively contraderived categories of co- respectively contramodules. We will show that curved Koszul triality is compatible with all of these binary operations.

\begin{definition}\label{binary functors}
Let $C$ be CDG coalgebra, $M$ a right CDG $C$-comodule, $N$ a left CDG $C$-comodule and $P$ a left CDG $C$-contramodule.
\begin{enumerate}
    \item The cotensor product $M\square_C N$ is the equalizer
    \begin{align*}
        M\otimes N\rightrightarrows M\otimes K\otimes N     
    \end{align*}
    of the right $C$-coaction on $M$ and the left $C$-coaction on $N$. The graded vector space $M\square_C N$ has a differential induced by the natural differential on $M\otimes N$, giving $M\square_C N$ the structure of a complex of vector spaces. We obtain a functor 
    \begin{displaymath}
        \square_C:comod\mhyphen C\times C\mhyphen comod\to Vect_k\ .
    \end{displaymath}
    
    \item The contratensor product $M\odot_C P$ is the coequalizer
    \begin{align*}
        M\otimes Hom(C,P)\rightrightarrows M\otimes P     
    \end{align*}
    of the right $C$-coaction on $M$ together with evaluation and the left $C$-contraaction on $P$. The graded vector space $M\odot_C P$ has a differential induced by the natural differential on $M\otimes P$, giving $M\odot_C P$ the structure of a complex of vector spaces. We obtain a functor 
    \begin{displaymath}
        \odot_C:comod\mhyphen C\times C\mhyphen contra\to Vect_k\ .
    \end{displaymath}
    
    \item The space of cohomomorphisms $Cohom_C(M,P)$ is the coequalizer
    \begin{align*}
        Hom(C\otimes N,P)\rightrightarrows Hom(N,P)     
    \end{align*}
    of the left $C$-coaction of $N$ and the left $C$-contraaction of $P$. The graded vector space $Cohom_C(M,P)$ has a differential induced by the natural differential on $Hom(N,P)$, giving $Cohom_C(M,P)$ the structure of a complex of vector spaces. We obtain a functor 
    \begin{displaymath}
        Cohom_C:comod\mhyphen C\times C\mhyphen contra\to Vect_k\ .
    \end{displaymath}
\end{enumerate}
\end{definition}

\begin{prop}[{\cite[Section 4.7.]{Positselski2011}}]\label{derived binary operations}
The functors defined in Definition~\ref{binary functors} induce derived functors
\begin{align*}
    \square^\mathbb{D}_C:D^{co}(comod\mhyphen C)\times D^{co}(C\mhyphen comod)\rightarrow D(Vect_k)\ ,\\
    \odot^\mathbb{D}_C:D^{co}(comod\mhyphen C)\times D^{ctr}(C\mhyphen contra)\rightarrow D(Vect_k)\ ,\\
    Cohom^\mathbb{D}_C:D^{co}(comod\mhyphen C)\times D^{ctr}(C\mhyphen contra)\rightarrow D(Vect_k)\ .
\end{align*}
\end{prop}

\begin{remark}
    In \cite[Section 4.7]{Positselski2011}, the derived functors of Proposition~\ref{derived binary operations} are referred to as $Cotor^C$, $Ctrtor^C$ and $Coext_C$, respectively.
\end{remark}

\begin{prop}\label{compatibility results}
    Let $A$ be a CDG algebra and $C$ be the extended bar construction of $A$. Then the following three diagrams commute:
    
    \begin{enumerate}
        \item$$\begin{tikzcd}
   D_c^\text{II}(mod\mhyphen A)\times D_c^\text{II}(A\mhyphen mod)\arrow[r, "F"]\arrow[d, "\otimes_A^{\text{II},\mathbb{L}}"]& D^{co}(comod\mhyphen C)\times D^{co}(C\mhyphen comod)\arrow[d, "\square_{C}^\mathbb{D}"]\\
   D(Vect_k)\arrow[r, "\cong"] & D(Vect_k)\ ,
\end{tikzcd}
$$
where $F=(-\otimes^\tau C,C\otimes^\tau-)$,

\item $$\begin{tikzcd}
   D_c^\text{II}(mod\mhyphen A)\times D_c^\text{II}(A\mhyphen mod)\arrow[r, "G"]\arrow[d, "\otimes_A^{\text{II},\mathbb{L}}"]& D^{co}(comod\mhyphen C)\times D^{ctr}(C\mhyphen contra)\arrow[d, "\odot_{C}^\mathbb{D}"]\\
   D(Vect_k)\arrow[r, "\cong"] & D(Vect_k)
\end{tikzcd}
$$
and

\item $$\begin{tikzcd}
   D_c^\text{II}(A\mhyphen mod)^{op}\times D_c^\text{II}(A\mhyphen mod)\arrow[r, "G"]\arrow[d, "\mathbb{R}Hom^{\text{II}}_A"]& D^{co}(C\mhyphen comod)^{op}\times D^{ctr}(C\mhyphen contra)\arrow[d, "Cohom_{C}^\mathbb{D}"]\\
   D(Vect_k)\arrow[r, "\cong"] & D(Vect_k)
\end{tikzcd}
$$
where $G=(-\otimes^\tau C,Hom^\tau(C,-))$.
    \end{enumerate}
\end{prop}
\begin{proof}
    We will show compatibility only with the derived cotensor product. The proofs of the two remaining assertions are analogous. Let $M$ be a right $A$-module and $N$ a left $A$-module then
    \begin{align*}
        \square^\mathbb{D}_{C}F(M,N)&=(M\otimes^\tau C)\square^\mathbb{D}_{C} (C\otimes^\tau N)\\
        &\cong M\otimes^\tau C\otimes^\tau N\\
        &\cong (M\otimes^\tau C\otimes^\tau A)\otimes_A (A\otimes^\tau C\otimes^\tau N)\\
        &=M\otimes_A^{\text{II},\mathbb{L}} N
    \end{align*}
    where the second isomorphism follows from the fact that $M\otimes^\tau C$ and $ C \otimes^\tau N$ are bifibrant objects in their respective model categories and the second to last isomorphism follows from the quasi-isomorphism $A\otimes^\tau C\simeq k$.
\end{proof}

\begin{remark}
    It was shown in \cite[Section 5.3]{Positselski2011} that comodule-contramodule correspondence turns the functor $\square_C$ into the functor $\odot_C$ and the functor $Cohom^\mathbb{D}_C$ into the functors $Ext_C:C\mhyphen comod\times C\mhyphen comod^{op}\to D(Vect_k)$ and $Ext^C:C\mhyphen contra\times C\mhyphen contra^{op}\to D(Vect_k)$. We therefore conclude that curved Koszul triality is compatible with all of the binary operations $\otimes^{\mathbb{L},\text{II}}_A$, $\square^\mathbb{D}_C$ and $\odot^\mathbb{D}_C$ as well as the corresponding derived hom-functors.
\end{remark}

\begin{remark}
    The main obstruction to promoting the results of Proposition~\ref{compatibility results} to statements about CDG bimodules stems from the fact that by Koszul duality there is an equivalence $D_c^\text{II}(A\otimes B^{op})\simeq D^{co}(\check{B}(A\otimes B^{op}))$, but $D^{co}(\check{B}(A\otimes B)\not\simeq D^{co}(\check{B}(A)\otimes \check{B}(B))$ in general. A construction addressing this problem is the subject of ongoing work.
\end{remark}

\Addresses


\begin{thebibliography}{99}

\bibitem{BalmerFavi2011}
P.~Balmer and G.~Favi, Generalized tensor idempotents and the telescope conjecture.
\emph{Proc. Lond. Math. Soc. (3)} \textbf{102} (2011), no.~6, 1161--1185


\bibitem{becker12}
H.~Becker, Models for singularity categories.
\emph{Adv. Math.} \textbf{254} (2014), 187--232


\bibitem{Beilinson96}
A.~Beilinson, V.~Ginzburg and W.~Soergel, Koszul duality patterns in representation theory.
\emph{J. Am. Math. Soc.} \textbf{9} (1996), no.~2, 473--527


\bibitem{Bondal_1991}
A.~I. Bondal and M.~M. Kapranov, Enhanced triangulated categories.
\emph{Math. USSR} \textbf{70} (1991), no.~1, 93--107


\bibitem{BoothLazarev2023_GlobalKoszulDuality}
M.~Booth and A.~Lazarev, Global Koszul duality.
2026, \url{https://arxiv.org/abs/2304.08409}

\bibitem{ChuangHolsteinLazarev2021}
J.~Chuan, J.~Holstein and A.~Lazarev, Maurer–Cartan moduli and theorems of Riemann–Hilbert type.
\emph{Appl. Categ. Struct.} \textbf{29} (2021), no.~4, 685--728


\bibitem{GuanHolsteinLazarev2025}
A.~Guan, J.~Holstein and A.~Lazarev, Hochschild cohomology of the second kind: Koszul duality and Morita invariance.
2026, \url{https://arxiv.org/abs/2312.16645}

\bibitem{GuanLazarev}
A.~Guan and A.~Lazarev, Koszul duality for compactly generated derived categories of second kind.
\emph{J. Noncommut. Geom.} \textbf{15} (2021), no.~4, 1355--1371


\bibitem{HessKedziorekRiehlShipley2017}
K.~Hess, M.~Kedziorek, E.~Riehl and B.~Shipley, A necessary and sufficient condition for induced model structures.
\emph{J. Topol.} \textbf{10} (2017), no.~2, 324--369


\bibitem{hinich}
V.~Hinich, Homological algebra of homotopy algebras.
\emph{Commun. Algebra} \textbf{25} (1997), no.~10, 3291--3323


\bibitem{Hovey1999}
M.~Hovey, \emph{Model categories}.
Math. Surv. Monogr. 63, American Mathematical Society (AMS), Providence, RI, 1999


\bibitem{Positselski_2022}
L.~Positselski, Contramodules.
\emph{Confluentes Math.} \textbf{13} (2021), no.~2, 93--182


\bibitem{Positselski2023KoszulDuality}
L.~Positselski, Differential graded Koszul duality: An introductory survey.
\emph{Bull. Lond. Math. Soc.} \textbf{55} (2023), no.~4, 1551--1640


\bibitem{Positselski_2010}
L.~Positselski, Homological algebra of semimodules and semicontramodules.
Monogr. Mat., Inst. Mat. PAN, N.S. 70, Birkhäuser, Cham, 2010


\bibitem{Positselski2011}
L.~Positselski, Two kinds of derived categories, Koszul duality, and comodule-contramodule correspondence.
\emph{Mem. Am. Math. Soc.} \textbf{212} (2011), no.~996


\bibitem{Positselski_2024}
L.~Positselski and J.~Šťovíček, Coderived and contraderived categories of locally presentable abelian DG-categories.
\emph{Math. Z.} \textbf{308} (2024), no.~14


\bibitem{ROOS_2006}
J.~E. Roos, Derived functors of inverse limits revisited.
\emph{J. Lond. Math. Soc., II.} \textbf{73} (2006), no.~1, 65--83


\bibitem{sweedler1969hopf}
M.~E. Sweedler, Hopf algebras.
Math. Lect. Note Ser., The Benjamin/Cummings Publishing Company, Reading, MA, 1969 


\end{thebibliography}
\end{document}